\documentclass[number,citesort,MSNbibl,dvips]{arxbj}
\usepackage{upgreek,multirow}
\usepackage{graphicx}

\aid{0}
\volume{18}
\issue{2}
\pubyear{2012}
\firstpage{391}
\lastpage{433}
\doi{10.3150/10-BEJ340}

\makeatletter
\newcommand{\fraca}[2]{{#1/#2}}
\newcommand{\fracb}[2]{{(#1)/#2}}
\newcommand{\fracc}[2]{{#1/(#2)}}
\newcommand{\fracd}[2]{{(#1/#2)}}
\newcommand{\frace}[2]{{(#1)/(#2)}}

\newtheorem{proposition}{Proposition}
\newtheorem{lemma}{Lemma}
\newtheorem{theorem}{Theorem}
\newtheorem{corollary}{Corollary}
\newremark{rem}{Remark}

\newcommand{\Charac}[1]{1_{\{#1\}}}
\newcommand{\cEstim}{\mathcal{E}}
\renewcommand{\epsilon}{\varepsilon}
\newcommand{\eps}{\epsilon}
\renewcommand{\phi}{\varphi}
\renewcommand{\tilde}{\widetilde}
 \newcommand{\imod}[1]{(\operatorname{mod}#1)}
\newcommand{\field}[1]{\mathbb{#1}}
\newcommand{\bR}{\field{R}}
\newcommand{\bS}{\field{S}}
\newcommand{\bN}{\field{N}}
\newcommand{\bE}{\field{E}}
\newcommand{\bL}{\field{L}}
\newcommand{\cA}{{\cal A}}
\newcommand{\cL}{{\cal L}}
\newcommand{\cM}{{\cal M}}
\newcommand{\cP}{{\cal P}}
\newcommand{\cV}{{\cal V}}
\newcommand{\cX}{{\mathcal X}}
\newcommand{\Arccos}{\operatorname{Arccos}}
\newcommand{\Card}{\operatorname{Card}}
\newcommand{\Diag}{\operatorname{Diag}}
\def\cM{{\field M}}
\def\cH{{\field H}}
\def\bR{{\field R}}
\def\a{\alpha}
\def\b{\beta}
\def\ONE{{ 1}}
\def\lambdaquadra{\omega}
\def\hf{\hat{f}}
\def\supp{\operatorname{supp}}

\newcommand{\eqref}[1]{(\ref{#1})}
\renewcommand{\emptyset}{\varnothing}
\newcommand{\cal}{\mathcal}
\makeatother

\begin{document}
\begin{frontmatter}

\title{Radon needlet thresholding}
\runtitle{Radon needlet thresholding}

\begin{aug}
\author{\fnms{G\'{e}rard} \snm{Kerkyacharian}\thanksref{e1}\ead[label=e1,mark]{kerk@math.jussieu.fr}},
\author{\fnms{Erwan} \snm{Le Pennec}\thanksref{e2}\ead[label=e2,mark]{lepennec@math.jussieu.fr}}
\and\\
\author{\fnms{Dominique} \snm{Picard}\corref{}\thanksref{e3}\ead[label=e3,mark]{picard@math.jussieu.fr}}
\runauthor{G. Kerkyacharian, E. Le Pennec and D. Picard}
\address{LPMA, CNRS, Universit{\'{e}} Paris Diderot, 175 rue du Chevaleret,
75013 Paris, France.\\
\printead{e1,e2,e3}}
\end{aug}

\received{\smonth{8} \syear{2009}}
\revised{\smonth{7} \syear{2010}}

%
\begin{abstract}
We provide a new algorithm for the treatment of the noisy inversion of
the Radon transform using an appropriate thresholding technique adapted
to a well-chosen new localized basis. We establish minimax results and
prove their optimality. In particular, we prove that the procedures
provided here are able to attain minimax bounds for any $\bL_p$ loss.
It is important to notice that most of the minimax bounds obtained here
are new to our knowledge. It is also important to emphasize the
adaptation properties of our procedures with respect to the regularity
(sparsity) of the object to recover and to inhomogeneous smoothness.
We perform a numerical study that is of importance since we especially
have to discuss the cubature problems and propose an averaging
procedure that is mostly in the spirit of the cycle spinning performed
for periodic signals.\looseness=-1
\end{abstract}

%

\end{frontmatter}

\section{Introduction}\label{sec1}
We consider the problem of inverting noisy observations of the
$d$-dimensional Radon transform. Obviously, the most immediate examples
occur for $d=2$ or 3. However, no major differences arise from
considering the general case.

There is considerable literature on the problem of reconstructing
structures from their Radon transforms, which is a fundamental problem in
medical imaging and, more generally, in tomography.
In our approach, we focus on several important points. We produce a
procedure that is efficient from an $\bL_2$ point of view, since this
loss function mimics quite well in many situations the preferences of
the human eye. On the other hand, we have at the same time the
requirement of clearly identifying the local bumps, of being able to
estimate the different level sets well. We also want the procedure to
enjoy good adaptation properties.
In addition, we require the procedure to be simple to implement.

At the heart of such a problem,
there is a notable conflict between the inversion part~-- which, in the
presence of noise, creates an instability reasonably handled by a
singular value decomposition (SVD) approach -- and
the fact that the SVD basis is very rarely localized and, as a
consequence, capable of representing local features of images,
which are especially important to recover.

Our strategy is to follow the approach started in \cite{5authors}, which
utilizes a construction borrowed from \cite{pxuball} (see also \cite
{pxukball})
of localized frames based on orthogonal polynomials on the ball, which
are closely related to
the Radon transform SVD basis.

To achieve the goals presented above, and especially adaptation to
different regularities and local inhomogeneous smoothness, we add a
fine-tuned subsequent thresholding process to the estimation performed in
\cite{5authors}.

This improves considerably the performances of the algorithm, both
from a theoretical point of view and a numerical point of view.
In effect, the new algorithm provides a~much better spatial adaptation,
as well as adaptation to the classes of regularity.
We prove here that the bounds obtained by the procedure are minimax
over a large class of Besov spaces and any $\bL_p$ loss: we provide
upper bounds for the performance of our algorithm and lower bounds for
the associated minimax rate.

It is important to notice that because we consider
different $\bL_p$ losses, we provide rates of convergence of new types
attained by our procedure. Those rates are minimax since they are
confirmed by
lower-bound inequalities.

The problem of choosing appropriate spaces of regularity on the ball
reflecting the standard objects analyzed in tomography is a highly
non-trivial problem. We decided to consider the spaces that seem to
stay closest to our natural intuition, that is, those that generalize
to the ball the approximation properties by polynomials.

The procedure gives very promising results in the simulation study. We
show that the estimates obtained by thresholding the needlets
outperform those obtained either by thresholding the SVD or by the
linear needlet estimate proposed in \cite{5authors}. An important
issue in the needlet scheme is the choice of the
quadrature in the needlet construction. We discuss the possibilities
proposed in the literature and consider a cubature formula based
on the full tensorial grid on the sphere, introducing an averaging
close to the cycle-spinning method.

Among others, one amazing result is that, to attain minimax
rates in the $\bL_\infty$ norm, we need to modify the
estimator. This result is also corroborated by the numerical results: see
Theorem \ref{mMinfty} and Figures \ref{fig:result} and \ref{fig:result2}.

In the first section, we introduce the Radon transform and the
associated SVD basis. The following section summarizes the construction
of the localized basis, the needlets.
The procedure is introduced in Section~\ref{sec4}, where the main theoretical
results are stated for
upper bounds and lower bounds.
Section~\ref{sec5} details the simulation study.
Section~\ref{sec6} details important properties of the needlet basis.
The proof of the two main
results stated in Section~\ref{sec4} is postponed in the two last sections.

\section{Radon transform and white noise model}\label{sec2}

\subsection{Radon transform}\label{sec2.1}

Here we recall the definitions and some basic facts about the Radon transform
(cf. \cite{Helgason,NATT,LOG}).
Denote by $B^d$ the unit ball in $\bR^d$,
that is,
$B^d = \{x =(x_1,\ldots,x_d)\in\bR^d\dvt |x| \le1\}$ with $|x|=(\sum
_{i=1}^d x_i^2)^{1/2}$
and, by $\bS^{d-1}$, the unit sphere in $\bR^d$.
The Lebesgue measure on~$B^d$ will be denoted by $\mathrm{d}x$
and the usual surface measure on $\bS^{d-1}$ by $\mathrm{d}\sigma(x)$
(sometimes we will also deal with the surface measure on $\bS^d$,
which will be denoted by $\mathrm{d}\sigma_d$).
We let~$|A|$ denote the measure $|A| = \int_A \mathrm{d}x$ if $A\subset B^d$
and $ |A| = \int_A \mathrm{d}\sigma(x) $ if $A\subset\bS^{d-1}$.

The Radon transform of a function $f$ is defined by
\[
Rf(\theta, s) =
\int_{\mathop{y\in\theta^\perp}\limits_{
s\theta+y\in B^d}} f(s\theta + y)\,\mathrm{d}y,
\qquad\theta\in\bS^{d-1},   s \in[-1,1],
\]
where $\mathrm{d}y$ is the Lebesgue measure of dimension $d-1$ and
$\theta^\perp=\{x \in\bR^d\dvt  \langle x , \theta\rangle= 0 \}$.
With a slight abuse of notation, we will rewrite this integral as
\[
Rf(\theta, s) = \int_{\langle y,\theta\rangle
=s}f(y)\,\mathrm{d}y.
\]
By Fubini's theorem, we have
\[
\int_{-1}^1 Rf (\theta, s)\,\mathrm{d}s = \int_{B^d} f(x)\,\mathrm{d}x.
\]

It is easy to see (cf. \cite{NATT}) that the Radon transform is a
bounded linear operator mapping
$\bL_2(B^d,\mathrm{d}x)$ into $\bL_2  (\bS^{d-1} \times [-1,1], \mathrm{d}\mu
(\theta,s) )$,
where
\[
\mathrm{d}\mu(\theta, s) = \mathrm{d}\sigma(\theta) \frac{\mathrm{d}s}{(1-s^2)^{(d-1)/2}}.
\]

\subsection{Noisy observation of the Radon transform}\label{sec2.2}

We consider observations of the form
\[
\mathrm{d}Y(\theta, s) = Rf (\theta, s)\,\mathrm{d}\mu(\theta, s) + \eps \,\mathrm{d}W(\theta, s),
\]
where the unknown function $f$ belongs to
$\bL_2(B^d,\mathrm{d}x) $.
The meaning of this equation is that, for any
$\phi(\theta,s)$ in $\bL_2(\bS^{d-1} \times [-1,1], \mathrm{d}\mu(\theta
,s) ),$
one can observe
\begin{eqnarray*}
Y_\phi
&=&\int\phi(\theta,s )\,\mathrm{d}Y(\theta,s) = \int_{\bS^{d-1} \times
[-1,1]} Rf(\theta,s)
\phi(\theta,s)\,\mathrm{d}\mu(\theta,s) + \eps\int\phi(\theta,s)\,\mathrm{d}W(\theta
,s)\\
&=& \langle Rf, \phi\rangle_\mu+ \eps W_\phi.
\end{eqnarray*}
Here $W_\phi=\int\phi(\theta,s)\,\mathrm{d}W(\theta,s)$ is a Gaussian field
of zero mean and covariance
\[
\bE(W_\phi,  W_\psi)= \int_{\bS^{d-1} \times [-1,1]} \phi(\theta
,s) \psi(\theta, s)\,\mathrm{d}\sigma(\theta)
\frac{\mathrm{d}s}{(1-s^2)^{(d-1)/2}} = \langle\phi, \psi\rangle_\mu.
\]
The goal is to recover the unknown function $f$ from the observation of $Y$.
Our idea is to refine the algorithms proposed in \cite{5authors}
using thresholding methods. In \cite{5authors}, estimation
schemes
are derived that combine
the stability and computability of SVD decompositions with
the localization and multiscale structure of wavelets.
To this end, a specific frame (essentially following the construction from
\cite{pxukball}) is used. It comprises elements of nearly exponential
localization and is, in addition, compatible with the SVD basis of the
Radon transform.

\subsection{Singular value decomposition of the Radon transform}\label
{SVD-Radon}

The SVD of the Radon transform was first established in
\cite{Davison,Louis}.
In this regard, we also refer the reader to \cite{NATT,xu7}.

\subsubsection{Jacobi and Gegenbauer polynomials}\label{Jacobi}\label{sec2.3.1}

The Radon SVD bases are defined in terms of Jacobi and Gegenbauer polynomials.
The Jacobi polynomials $P_n^{(\alpha,\beta)}$, $n\ge0$,
constitute an orthogonal basis for the space
$\bL_2([-1, 1], w_{\alpha,\beta}(t)\,\mathrm{d}t)$ with weight $w_{\alpha
,\beta}(t)=(1-t)^\a(1+t)^\b$,
$\a, \b>-1$.
They are standardly normalized by
$P_n^{(\a,\b)}(1)={{n+\a\choose n}}$ and then \cite{AAR,Erdelyi,SZG}
\[
\int_{-1}^1 P_n^{(\alpha,\beta)}(t)
P_{n'}^{(\a,\b)}(t)w_{\alpha,\beta}(t)\,\mathrm{d}t = \delta_{n, n'}
h_n^{(\alpha,\beta)},
\]
where
\[
h_n^{(\alpha,\beta)} =
\frac{2^{\alpha+\beta+1}}{(2n+\alpha+\beta+1)}
\frac{\Gamma(n+\a+1)\Gamma(n+\b+1)}{\Gamma(n+1)\Gamma(n+\a+\b+1)}.
\]
The Gegenbauer polynomials $C_n^\lambda$ are a particular case of
Jacobi polynomials
and are traditionally defined by
\[
C_n^\lambda(t)
=\frac{(2\lambda)_n}{(\lambda+1/2)_n}
P_n^{(\lambda-1/2, \lambda-1/2)}(t),
\qquad\lambda>-1/2,
\]
where, by definition, $(a)_n= a(a+1)\cdots (a+n-1)=\frac{\Gamma
(a+n)}{\Gamma(a)}$
(note that in \cite{SZG} the Gegenbauer polynomial $C_n^\lambda$ is
denoted by $P_n^\lambda$).
It is readily seen that
$C_n^\lambda(1)= {n+2\lambda-1 \choose n}=\frac{\Gamma(n+2\lambda
)}{n!\Gamma(2\lambda)}$ and
\[
\int_{-1}^1 C^\lambda_n (t) C^\lambda_m(t) (1-t^2)^{\lambda- 1/2}\,\mathrm{d}t
=\delta_{n, m}h_n^{(\lambda)}
\qquad\mbox{with }
h_n^{(\lambda)}= \frac{2^{1-2\lambda}\pi}{\Gamma(\lambda)^2}
\frac{\Gamma(n+2\lambda)}{(n+\lambda)\Gamma(n+1)}.
\]

\subsubsection{Polynomials on   $B^d$ and $\bS^{d-1}\label{sec2.3.2}
$}
\label{polynom}

Let $\Pi_n(\bR^d) $ be the space of all polynomials in $d$ variables
of degree $\le n$.
We denote by~$\cP_n(\bR^d)$ the space of all homogeneous polynomials
of degree $n$
and by $\cV_n(\bR^d)$ the space of all polynomials of degree $n$
that are orthogonal to lower-degree polynomials with respect to
the
Lebesgue measure on $B^d$.
Of course, $\cV_0(\bR^d)$ will be the set of all constants.
We have the following orthogonal decomposition:
\[
\Pi_n(\bR^d) =\bigoplus_{k=0}^n \cV_k(\bR^d).
\]

Also, denote by $\cH_n(\bR^d)$ the subspace of all harmonic
homogeneous polynomials of degree $n$
(i.e., $ Q \in\cH_n(\bR^d)$ if $Q\in\cP_n(\bR^d)$ and $\Delta
Q =0$)
and, by $\cH_n(\bS^{d-1}),$ the (injective) restriction of the
polynomials from $\cH_n(\bR^d)$ to $\bS^{d-1}$.
It is well known that
\[
N_{d-1}(n)
=\dim(\cH_n(\bS^{d-1}))
= \pmatrix{n+d-1\cr d-1} - \pmatrix{n+d-3 \cr d-1}
\sim n^{d-2}.
\]
Let $\Pi_n(\bS^{d-1}) $ be the space of restrictions to $\bS^{d-1}$
of polynomials of degree
$\leq n$ on $\bR^d$.
It is also well known that
\[
\Pi_n(\bS^{d-1})= \bigoplus_{m=0}^n \cH_{m}(\bS^{d-1})
\]
(the orthogonality is, with respect to the surface measure, $\mathrm{d}\sigma$
on $\bS^{d-1}$).
$\cH_{l}(\bS^{d-1})$ is called the space of spherical harmonics of
degree $d$ on the sphere
$\bS^{d-1}$.

Let $Y_{l,m}$, $1\leq m \leq N_{d-1}(l)$, be an orthonormal basis of
$\cH_{l}(\bS^{d-1})$, that is,
\[
\int_{\bS^{d-1}}Y_{l,m}(\xi) \overline{Y_{l,m'}(\xi)}\,\mathrm{d}\sigma(\xi
) = \delta_{m,m'}.
\]
Then the natural extensions of $Y_{l,m}$ on $B^d$ are defined by
$Y_{l,m}(x)=|x|^l Y_{l,m} ( \frac{x}{|x|} )$ and satisfy
\begin{eqnarray*}
\int_{B^{d}}Y_{l,m}(x) \overline{Y_{l,m'}(x)}\,\mathrm{d}x
&=& \int_0^{1} r^{d-1} \int_{\bS^{d-1}}Y_{l,m}(r\xi) \overline
{Y_{l,m'}(r\xi)}\,\mathrm{d}\sigma(\xi)\,\mathrm{d}r \\
&=& \int_0^{1} r^{d+2l-1} \int_{\bS^{d-1}}Y_{l,m}(\xi) \overline
{Y_{l,m'}(\xi)}\,\mathrm{d}\sigma(\xi)\,\mathrm{d}r
= \delta_{m,m'} \frac1{2l+d}.
\end{eqnarray*}
For more details, we refer the reader to \cite{DUXU}.

The spherical harmonics on $\bS^{d-1}$ and orthogonal polynomials on
$B^d$ are naturally related
to Gegenbauer polynomials.
Thus, the kernel of the orthogonal projector onto $\cH_n(\bS^{d-1})$
can be written as (see \cite{STW}):
%
\begin{equation}
\label{legen}
\sum_{m=1}^{N_{d-1}(n)} Y_{l, m}(\xi) \overline{Y_{l, m}(\theta)}
= \frac{2n+ d-2}{(d-2)|\bS^{d-1}|} C^{\fracb{d-2}{2}}_n (\langle\xi
, \theta\rangle).
\end{equation}
The ``ridge'' Gegenbauer polynomials $C_n^{d/2}(\langle x, \xi \rangle
)$ are orthogonal
to $\Pi_{n-1}(B^d)$ in $\bL_2(B^d)$ and the kernel $L_n(x,y)$ of the
orthogonal projector onto $\cV_n(B^d)$
can be written in the form (see \cite{petrush,xu7})
%
\begin{eqnarray}
\label{orth-projector}
L_n(x,y)
&=& \frac{2n+d}{|\bS^{d-1}|^2}
\int_{\bS^{d-1}}C^{d/2}_n ( \langle x, \xi \rangle)C^{d/2}_n (
\langle y, \xi \rangle)\,\mathrm{d}\sigma(\xi)\nonumber
\\[-8pt]
\\[-8pt]
&=& \frac{(n+1)_{d-1}}{2^d \uppi^{d-1}}
\int_{\bS^{d-1}}\frac{C^{d/2}_n ( \langle x, \xi \rangle)
C^{d/2}_n ( \langle y, \xi \rangle) }{\| C^{d/2}_n\|^2}\,\mathrm{d}\sigma(\xi
).
\nonumber
\end{eqnarray}

The following important identities are valid for ``ridge'' Gegenbauer
polynomials:
%
\begin{equation}
\label{ridge-Gegen1}
\int_{B^d} C^{d/2}_n (\langle\xi, x \rangle)C^{d/2}_n (\langle\eta
, x \rangle)\,\mathrm{d}x
=\frac{h_n^{(d/2)}}{C_n^{d/2}(1)}C^{d/2}_n (\langle\xi, \eta\rangle),
\qquad\xi, \eta\in\bS^{d-1},
\end{equation}
and, for $x\in B^d$, $\eta\in\bS^{d-1}$,
%
\begin{equation}
\label{ridge-Gegen2}
\int_{\bS^{d-1}} C^{d/2}_n (\langle\xi, x \rangle)C^{d/2}_n
(\langle\xi, \eta\rangle)\,\mathrm{d}\sigma(\xi)
=|\bS^{d-1}| C^{d/2}_n (\langle\eta, x \rangle);
\end{equation}
see \cite{petrush}.
By (\ref{orth-projector}) and (\ref{ridge-Gegen2})
\[
L_n(x,\xi)
= \frac{(2n+d)}{|\bS^{d-1}|} C^{d/2}_n ( \langle x, \xi \rangle),
\qquad\xi\in\bS^{d-1},
\]
and again by (\ref{orth-projector})
\[
\int_{\bS^{d-1}}L_n (x,\xi) L_n(y, \xi)\,\mathrm{d}\sigma(\xi) = (2n+d)L_n(x,y).
\]
%
\subsubsection{The SVD of the Radon transform}\label{SVD}\label{sec2.3.3}

Assume that
$\{Y_{l,m}\dvt 1\leq m \leq N_{d-1}(l)\}$ is an orthonormal basis for $\cH
_{l}(\bS^{d-1})$.
Then it is standard and easy to see that the family of polynomials,
\begin{eqnarray*}
&&f_{k,l,m} (x)
= (2k+d)^{1/2}P_j^{(0,  l + d/2 -1)} (2|x|^2-1)Y_{l, m}(x),
 \\
 && \quad    0\leq l \leq k,   k-l =2j,    1\leq i \leq N_{d-1}(l),
\end{eqnarray*}
form an orthonormal basis of $\cV_k(B^d)$; see \cite{DUXU}.
Here, as before, $Y_{l, m}(x)=|x|^lY_{l, m}(x/|x|)$.
On the other hand, the collection
\[
g_{k,l,m}(\theta, s)
= \bigl[h_k^{(d/2)}\bigr]^{-1/2}(1-s^2)^{(d-1)/2} C^{d/2}_k(s) Y_{l,m }(\theta),
\qquad k\ge0,    l \ge0,    1\leq m \leq N_{d-1}(l),
\]
is an orthonormal basis of $\bL_2(\bS^{d-1}\times[-1,1], \mathrm{d}\mu
(\theta,s))$.
Most important, the Radon transform
$R\dvtx \bL_2(B^d) \mapsto\bL_2(\bS^{d-1}\times[-1,1], \mathrm{d}\mu(\theta,s))$
is a one-to-one mapping and
\[
Rf_{k,l,m} =\lambda_k g_{k,l,m}, \qquad
R^*g_{k,l,m} =\lambda_k f_{k,l,m},
\]
where
\[
\lambda_k^2 = \frac{2^d\pi^{d-1}}{(k+1)(k+2)\cdots (k+d-1)}
= \frac{2^d\pi^{d-1}}{(k+1)_{d-1}}\sim k^{-d+1}.
\]
More precisely, we have
for any $f\in\bL_2(B^d)$
\[
Rf=
\sum_{k\ge0}  \lambda_k
\sum_{0 \leq l \leq k,  k-l \equiv0\ \imod{2}}
\sum_{ 1 \leq m \leq N_{d-1}(l) }
\langle f, f_{k, l, m}\rangle g_{k,l,m}.
\]
Furthermore, for $f\in\bL_2(B^d),$
\[
f= \sum_{k\ge0}  \lambda_k^{-1}
\sum_{0 \leq l \leq k,  k-l \equiv0\ \imod{2}}
\sum_{ 1 \leq m \leq N_{d-1}(l) }
\langle Rf, g_{k, l, m}\rangle_\mu f_{k,l,m}.
\]
In the above identities, the convergence is in $\bL_2$.

For the Radon SVD, we refer the reader to
\cite{NATT,xu7,5authors}.

\section{Construction of needlets on the ball}\label{sec3}

In this section, we briefly recall the construction of the needlets on
the ball.
This construction is due to \cite{pxuball}. Its aim is to build a
very well-localized tight frame constructed using the eigenvectors of
the Radon transform. For more precision, we refer the reader to \cite
{pxuball,kyoto,5authors}

Let $\{f_{k,l,m}\}$ be
the orthonormal basis of $\cV_k(B^d)$ defined in Section \ref{SVD}.
Denote by $T_k$ the index set of this basis, that is,
$
T_k=\{(l, m)\dvt  0\leq l \leq k,  l\equiv k\ \imod{2},   0\leq m \leq
N_{d-1}(l)\}$.
Then the orthogonal projector of $\bL_2(B^d)$ onto $\cV_k(B^d)$
can be written in the form
\[
L_kf = \int_{B^d} f(y) L_k(x,y)\,\mathrm{d}y
\qquad\mbox{with }
L_k(x,y) = \sum_{l,m \in T_k} f_{k,l,m}(x) f_{k,l,m}(y).
\]
Using (\ref{legen}), $L_k(x,y)$ can be written in the form
\begin{eqnarray*}
&&L_k(x,y)\\
&& \quad  =(2k+d) \sum_{l\le k,  k-l\equiv0\ \imod{2}}
P_j^{(0, l +d/2-1)} (2|x|^2-1) |x|^l
P_j^{(0, l +d/2-1)} (2|y|^2-1)|y|^l \\
&&\hphantom{=(2k+d) \sum_{l\le k,  k-l\equiv0\ \imod{2}}}   \quad {}  \times
\sum_{m} Y_{l, m} \biggl(\frac x{|x|} \biggr)Y_{l,m} \biggl(\frac y{|y|}
\biggr) \\
&& \quad = \frac{ (2k+d)}{|\bS^{d-1}|} \sum_{l\le k,  k-l\equiv0\ \imod
{2}}
P_j^{(0, l +d/2-1)} (2|x|^2-1)|x|^l P_j^{(0, l +d/2-1)} (2|y|^2-1)
|y|^l \\
&& \hphantom{= \frac{ (2k+d)}{|\bS^{d-1}|} \sum_{l\le k,  k-l\equiv0\ \imod
{2}}}\quad {} \times \biggl(1+\frac l{d/2-1} \biggr)
C^{d/2-1}_l  \biggl( \biggl\langle\frac x{|x|}, \frac y{|y|}  \biggr\rangle
 \biggr).
\end{eqnarray*}
Another representation of $L_k(x,y)$ has already been given in (\ref
{orth-projector}).
Clearly,
%
\begin{equation}
\label{pro}
\int_{B^d} L_k(x,z) L_{k'}(z,y)\,\mathrm{d}z = \delta_{k,k'} L_k(x,y)
\end{equation}
and, for $f\in\bL_2(B^d),$
%
\begin{equation}
\label{orthog-dec}
f=\sum_{k\ge0} L_kf
\quad\mbox{and}\quad
\|f\|_2^2 =\sum_{k} \|L_k f\|_2^2 =\sum_{k} \langle L_k f, f \rangle.
\end{equation}

The construction of the needlets is based on the classical
Littlewood--Paley decomposition and a subsequent
discretization.

Let $a\in C^\infty[0, \infty)$ be a cut-off function such that
$0 \leq a \leq1$,
$a(t)=1$ for $t\in[0, 1/2]$ and $\supp a \subset[0, 1]$.
We next use this function to introduce a sequence of operators on $\bL_2(B^d)$.
For $j\ge0$, write
\begin{eqnarray*}
&&A_j f (x) = \sum_{k\ge0} a \biggl(\frac k{2^j} \biggr) L_k f(x) = \int
_{B^d} A_j(x,y)f(y)\,\mathrm{d}y\\
&& \quad \mbox{with }   A_j(x,y)= \sum_{k} a \biggl(\frac k{2^j} \biggr)
L_k (x,y).
\end{eqnarray*}
Also, we define $B_jf = A_{j+1}f-A_jf$.
Then, setting $b(t) = a(t/2) - a(t),$ we have
\begin{eqnarray*}
&&B_j f (x) = \sum_{k} b \biggl(\frac k{2^j} \biggr)L_k f(x) = \int_{B^d}
B_j(x,y)f(y)\,\mathrm{d}y\\
&& \quad \mbox{with }
B_j(x,y)= \sum_{k} b \biggl(\frac k{2^j} \biggr)L_k (x,y).
\end{eqnarray*}
Obviously, for $f\in\bL_2(B^d)$,
\[
\langle A_jf,f \rangle= \sum_{k} a \biggl(\frac k{2^j} \biggr) \langle
L_k f, f \rangle\leq \| f\|_2^2.
\]

An important result from \cite{pxuball} (cf. \cite{pxukball})
asserts that
the kernels $A_j(x, y)$, $B_j(x, y)$ have nearly exponential
localization. Namely,
for any $M>0$ there exists a constant $C_M>0$ such that
%
\begin{equation}\label{EqFond}
|A_j(x,y)|, |B_j(x,y)| \leq C_M \frac{ 2^{jd}}{(1+2^jd(x,y))^M \sqrt{
W_j(x)} \sqrt{ W_j(y)} },
\qquad x, y\in B^d,
\end{equation}
where
%
\begin{equation}
W_j(x) = 2^{-j} +\sqrt{ 1-|x|^2}, \qquad|x|^2 =|x|^2_{d} =\sum
_{i=1}^d x_i^2
\label{Wj}
\end{equation}
and
\[
d(x,y)= \Arccos\bigl( \langle x,y\rangle+ \sqrt{ 1-|x|^2}\sqrt{ 1-|y|^2}\bigr),
\qquad
\langle x,y\rangle= \sum_{i=1}^d x_iy_i.
\]

Let us define
\[
C_j(x,y) =\sum_{k} \sqrt{a \biggl(\frac k{2^j} \biggr)} L_k (x,z)
\quad\mbox{and}\quad
D_j(x,y)=\sum_{k} \sqrt{b \biggl(\frac k{2^j} \biggr)} L_k (x,z).
\]
Note that $C_j$ and $D_j$ have the same localization as the
localization of $A_j$, $B_j$
in (\ref{EqFond}) (cf.~\cite{pxuball}).
Using (\ref{pro}), we get,
%
\begin{equation}
\label{D1}
A_j(x,y) = \int_{B^d} C_j(x,z) C_j(z,y)\,\mathrm{d}z
,\qquad
B_j(x,y) = \int_{B^d} D_j(x,z) D_j(z,y)\,\mathrm{d}z.
\end{equation}
And, obviously,
$z \mapsto C_j(x,z) C_j(z,y)$ (resp., $D_j(x,z) D_j(z,y)$)
are polynomial of degrees $< 2^{j+1}$.

The following proposition follows from results in \cite{pxuball} and
\cite{cubxu} and establishes a cubature formula.

\begin{proposition}\label{prop:CUB}
Let $\{B(\tilde\xi_i,\rho)\dvt  i \in I\}$ be a maximal family of
disjoint spherical caps of radius $\rho= \tau2^{-j}$ with
centers on the hemisphere $\bS_+^d$.
Then for sufficiently small $0<\tau\le1$ the set of points
$\chi_j=\{\xi_i\dvt  i \in I \}$
obtained by projecting the set $\{\tilde\xi\dvt  i \in I\}$ on $B^d$
is a set of nodes of a cubature formula that is exact for
$\Pi_{2^{j+2}}(B^d)$:
for any $P\in \Pi_{2^{j+2}}(B^d)$,
\[
\int_{B^d} P(u)\,\mathrm{d}u =
\sum_{\xi \in \chi_j} \lambdaquadra_{j,\xi} P(\xi),
\]
where,
moreover, the coefficients $\lambdaquadra_{j,\xi}$ of this cubature
are positive and
satisfy $\lambdaquadra_{j,\xi}\sim W_j(\xi)2^{-jd}$, and the
cardinality of the set $\chi_j$ is of order $ 2^{jd}$.
\end{proposition}

\subsection{Needlets}\label{needlets}\label{sec3.0.4}

Going back to identities \eqref{D1} and applying the cubature formula
described in Proposition~\ref{prop:CUB}, we get
\begin{eqnarray*}
A_j(x,y) &=& \int_{B^d} C_j(x,z) C_j(z,y)\,\mathrm{d}z
=\sum_{\xi\in\chi_j} \lambdaquadra_{j,\xi} C_j (x,\xi)C_j (
y,\xi)
\quad\mbox{and}\\
B_j(x,y) &=& \int_{B^d} D_j(x,z) D_j(z,y)\,\mathrm{d}z
=\sum_{\xi\in\chi_j} \lambdaquadra_{j,\xi} D_j (x,\xi)D_j (
y,\xi).
\end{eqnarray*}
We define the \textit{father needlets}, $\phi_{j,\xi}$, and
the \textit{mother needlets}, $\psi_{j,\xi}$, by
\[
\phi_{j,\xi} (x) = \sqrt{ \lambdaquadra_{j,\xi}} C_j (x,\xi)
\quad\mbox{and}\quad
\psi_{j,\xi} (x) = \sqrt{ \lambdaquadra_{j,\xi}} D_j (x,\xi),
\qquad\xi\in\chi_j,   j\ge0.
\]
We also set
$\psi_{-1,0}=\frac{\ONE_{B^d}}{|B^d|}$ and $\chi_{-1}=\{0\}$.
From above, it follows that
\[
A_j(x,y) =\sum_{\xi\in\chi_j} \phi_{j,\xi} (x) \phi_{j,\xi} (y),
\qquad
B_j(x,y) =\sum_{\xi\in\chi_j} \psi_{j,\xi} (x) \psi_{j,\xi} (y).
\]
Therefore,
\[
A_j f (x) = \int_{B^d} A_j(x,y) f(y)\,\mathrm{d}y=
\sum_{\xi\in\chi_j} \langle f, \phi_{j,\xi} \rangle \phi_{j,\xi}
=\sum_{\xi\in\chi_j} \alpha_{j,\xi} \phi_{j,\xi},
\qquad\alpha_{j,\xi}= \langle f, \phi_{j,\xi} \rangle
\]
and
\[
B_j f (x) = \int_{B^d} B_j(x,y) f(y)\,\mathrm{d}y=
\sum_{\xi\in\chi_j} \langle f, \psi_{j,\xi} \rangle \psi_{j,\xi}
=\sum_{\xi\in\chi_j} \beta_{j,\xi} \psi_{j,\xi},
\qquad\beta_{j,\xi}=\langle f, \psi_{j,\xi} \rangle.
\]

It is easy to prove (see \cite{pxuball}) that
\[
\|\phi_{j,\xi} \|_2 \leq1.
\]
From (\ref{orthog-dec}) and the fact that
$\sum_{j\ge0} b(t2^{-j})=1$ for $t\in[1, \infty)$, it readily
follows that
\[
f = \sum_{j\ge-1} \sum_{\xi\in
\chi_j } \langle f, \psi_{j,\xi} \rangle \psi_{j,\xi},
\qquad f\in\bL_2(B^d),
\]
and, taking the inner product with $f$, it leads to
\[
\|f\|_2^2 = \sum_j \sum_{\xi\in\chi_j } |\langle f, \psi_{j,\xi}
\rangle|^2.
\]
In turn, this shows that the family $\{\psi_{j,\xi}\}$ is a tight
frame for $\bL_2(B^d)$.

\section{Needlet inversion of a noisy Radon transform and minimax performances}\label{sec4}

Our estimator is based on an appropriate thresholding of a needlet
expansion as follows.
$f$ can be decomposed using the frame above:
\[
f = \sum_{j\ge-1} \sum_{\xi\in
\chi_j } \langle f, \psi_{j,\xi} \rangle \psi_{j,\xi}.
\]

Our estimation procedure will be defined by the following steps
%
\begin{eqnarray}
\label{step1}
 \widehat{\alpha}_{k,l,m} &=& \frac1{\lambda_k }\int
g_{k,l,m}\,\mathrm{d}Y,
\\
\label{estimator-alpha}
 \widehat{\beta}_{j, \xi}
&=& \sum_{k, l, m}
\gamma^{j,\xi}_{k,l,m} \widehat{\alpha}_{k,l,m}
\end{eqnarray}
with
\[
\gamma^{j,\xi}_{k,l,m}=\langle g_{k,l,m},\psi_{j,\xi
}\rangle
\]
 and
\begin{equation}\label{estim}
 \hat f=\sum_{j=-1}^{J_\eps} \sum_{\xi\in
\chi_j }\widehat{\beta}_{j, \xi}\ONE_{\{|\widehat{\beta}_{j, \xi
}|\ge\kappa2^{j\nu}c_\eps\}}\psi_{j,\xi}
\end{equation}
with
%
\begin{equation}\nu=(d-1)/2.\label{nu}
\end{equation}
Hence, our procedure has three steps: the first one (\ref{step1})
corresponds to the inversion of the operator in the SVD basis, the
second one (\ref{estimator-alpha}) projects on the needlet basis and
the third one (\ref{estim}) ends up the procedure with a final thresholding.
The tuning parameters of this estimator are
\begin{itemize}
\item The range $J_\eps$ of resolution levels will be taken such that
\[
2^{J_\eps
(d-\fraca12)} \leq  \bigl(\eps\sqrt{\log1/\eps} \bigr)^{-1} <
2^{(J_\eps+1)
(d-\fraca12)}.
\]
\item The threshold constant $\kappa$ is an important tuning of our
method. The theoretical point of view asserts that, for $\kappa$ above
a constant (for which our evaluation is probably not optimal), the
minimax properties hold. Evaluations of $\kappa$ from the simulation
point of view are also given.
\item$c_\eps$ is a constant depending on the noise level. We shall
see that the following choice is appropriate:
\[
c_\eps=\eps\sqrt{\log1/\eps}.
\]
\item Notice that the threshold function for each coefficient contains
$2^{j\nu}$. This is due to the inversion of the Radon operator and the
concentration relative to the $g_{k,l,m}$'s of the needlets.
\item It is important to remark here that, unlike the (linear)
procedures proposed in \cite{5authors}, this one does not require the
knowledge of the regularity while, as will be seen in the sequel, it
attains bounds that are as good as the linear ones and even better
since they are handling much wider ranges for the parameters of the
Besov spaces.
\end{itemize}

We will consider the minimax properties of this estimator on the Besov
bodies constructed on the needlet basis. In \cite{pxukball}, it is
proved that
these spaces can also be described as approximation spaces, so they
have a genuine meaning
and can be compared to standard Sobolev spaces.

We define here the Besov body, $B^s_{\pi,r}$, as the space of functions
$f=\sum_{j\ge-1} \sum_{\xi\in
\chi_j } \beta_{j,\xi} \psi_{j,\xi}$
such that
\[
\sum_j2^{jsr} \biggl(\sum_{\xi\in\chi_j} (|\beta_{j,\xi}| \|\psi
_{j,\xi}\|_\pi)^\pi \biggr)^{r/\pi} <\infty
\]
(with the obvious modifications for the cases $\pi$ or $r=\infty$)
and $B^s_{\pi,r}(M),$ the ball of radius~$M$ of this space.

\begin{theorem} \label{mM}
For $0< r\le\infty$, $\pi\ge1$, $1\le p< \infty$, there exist some
constant $c_p=c_p(s,\pi,r,M)$,  $\kappa_0$ such that if
$\kappa\ge\kappa_0$, $s>(d+1)(\frac{1}\pi-\frac1p)_+$ and, in
addition, if $\pi< p$, $s>\frac{d+1}\pi-\frac12$:

\begin{longlist}[(3)]
\item[(1)]
If $\frac1p <\frac d{d+1}$,
\begin{eqnarray*}
&&\sup_{f\in B^s_{\pi, r}(M)} (\bE\|\hf-f\|_p^p )^{\fraca1p}\\
&& \quad
\leq c_p  ({\log1/\epsilon} )^{\fraca{p}{2}}\\
&& \qquad {}\times\bigl ( \epsilon\sqrt{\log1/\epsilon} \bigr)^{\frace{s-
(d+1)(1/\pi-1/p)}{s +d -(d+1)/\pi}\wedge
\fracc s{s+ d-1/2}\wedge\frace{s-2(1/\pi-1/p)}{s+d-2/\pi}}.
\end{eqnarray*}
\item[(2)]
If $\frac d{d+1} \leq\frac{1}{p}$
and $d>2$ or $p>1$,
\[
\sup_{f\in B^s_{\pi, r}(M)} (\bE\|\hf-f\|_p^p )^{\fraca1p}
\leq c_p  ({\log1/\epsilon} )^{\fraca{p}{2}}
\bigl ( \epsilon\sqrt{\log1/\epsilon} \bigr)^{
\fracc s{s+ d-1/2}\wedge\frace{s-2(1/\pi-1/p)}{s+d-2/\pi}}.
\]
\item[(3)]
If $d=2$ and $p=1,$
\[
\sup_{f\in B^s_{\pi, r}(M)} (\bE\|\hf-f\|_1 )
\leq c_1  ({\log1/\epsilon} )^{\fraca{1}{2}}  \bigl(
\epsilon\sqrt{\log1/\epsilon} \bigr)^{\fracc s{s+ 2-1/2}}.
\]
\end{longlist}
\end{theorem}

\begin{figure}
\includegraphics{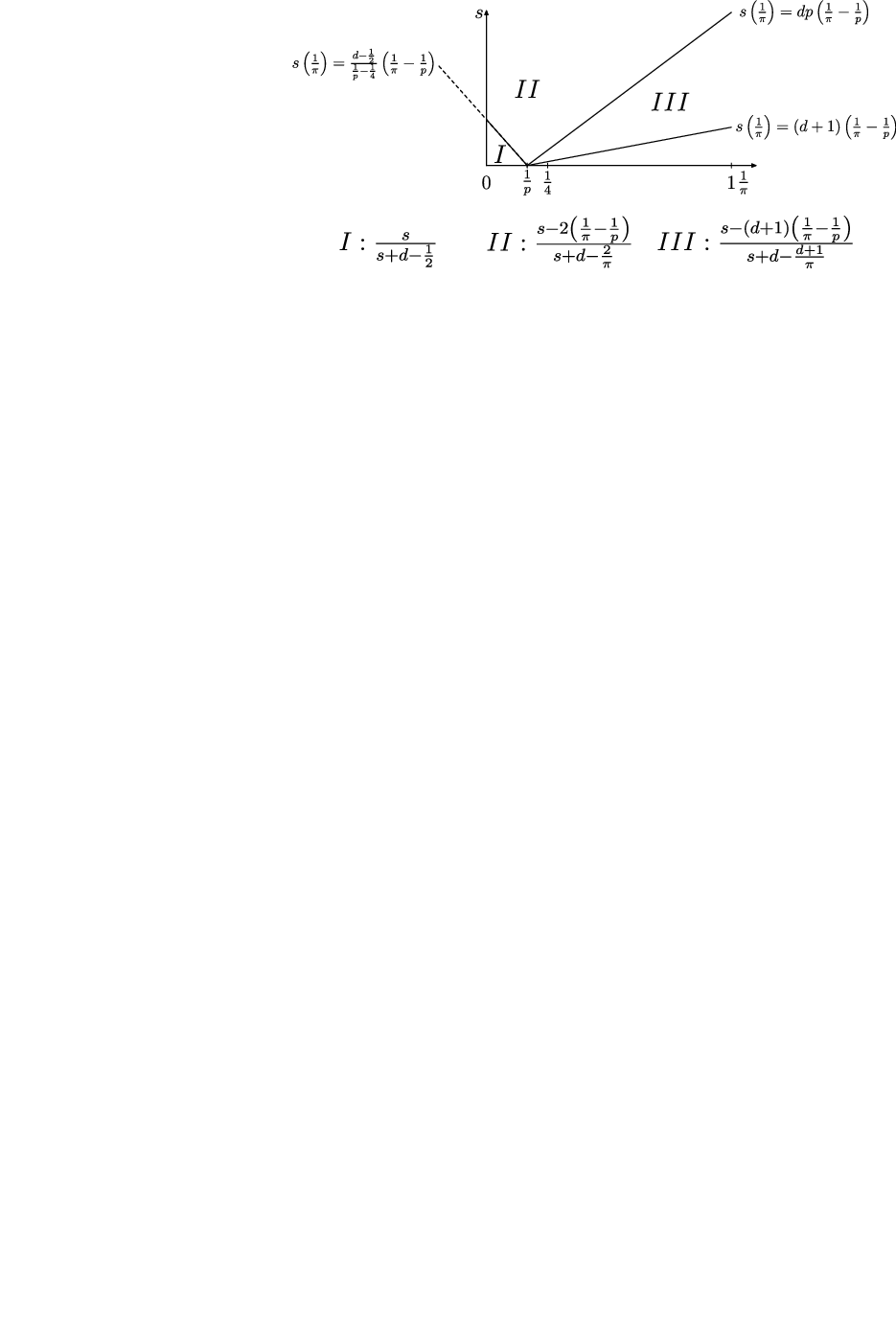}%
\caption{The three different minimax rate type zones are shown with
respect to the Besov space\vspace*{1pt} parameters
$s$ and $\pi$ for a fixed loss norm $L^p$ with $0 < \frac{1}{p} <
\frac{1}{4}$.}
\label{fig:lowerbound1}
\end{figure}

\begin{figure}[b]

\includegraphics{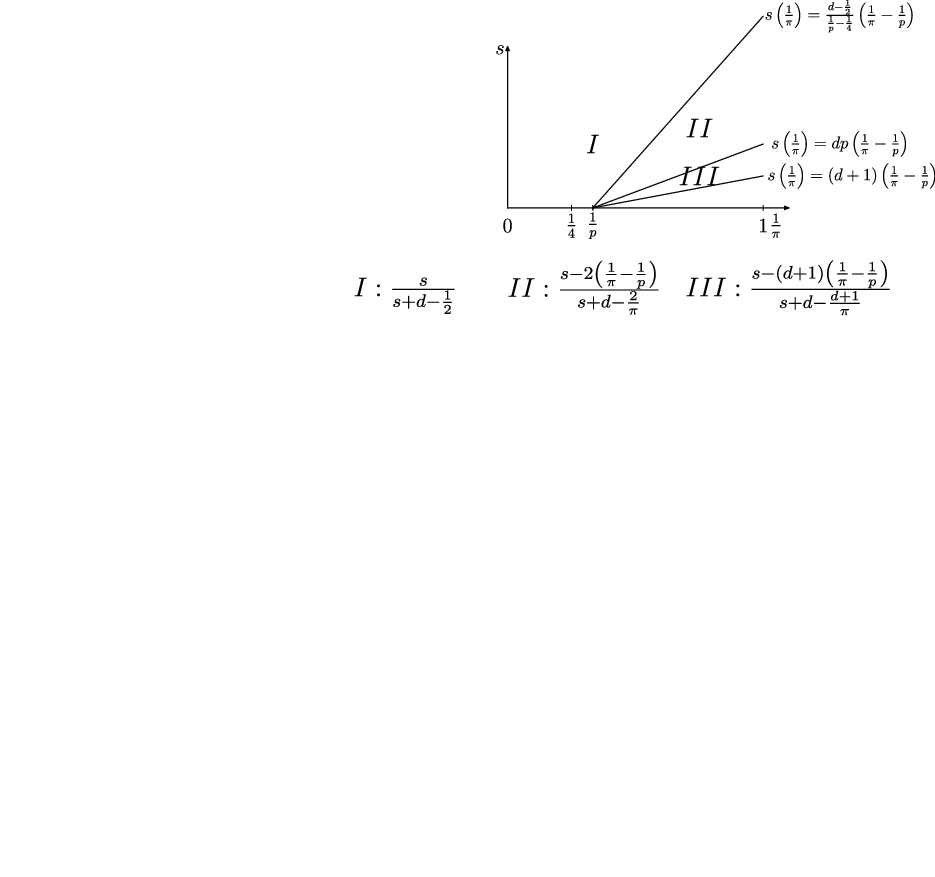}%
\caption{The three different minimax rate type zones are shown with respect to
the Besov space parameters
$s$ and $\pi$ for a fixed loss norm $p$ with $\frac{1}{4} < \frac
{1}{p} <
\frac{{d}}{d+1}$.}
\label{fig:lowerbound2}
\end{figure}

\begin{rem}
Up to logarithmic terms, the rates observed here are minimax, as will
appear in the following theorem. It is
known that in this kind of estimation, full adaptation yields
unavoidable extra logarithmic terms. The rates of the logarithmic terms
obtained in these
theorems are, most of the time, suboptimal (for instance, for obvious
reasons, the case $p=2$ yields fewer logarithmic terms). A more
detailed study could lead to optimized rates, which we decided not to
include here for the sake of simplicity.

The cumbersome comparisons of the different rates of convergence are
summarized in Figures~\ref{fig:lowerbound1} and~\ref{fig:lowerbound2}
for the case $0< \frac1p <\frac d{d+1}$.
These figures illustrate and highlight the differences between the
cases $p>4$ and $p<4$. We put $\frac1p$ as the horizontal axis and the
regularity $s$ as the vertical axis.
As explained later, after the lower-bound results, zones~I and II
correspond to two different types of the so called ``dense'' case,
whereas zone III corresponds to the ``sparse'' case.
\end{rem}
%

For the case of an $\bL_\infty$ loss function, we have a slightly
different result since the threshol\-ding depends on the $\bL_\infty$
norm of the local needlet.
Let us consider the following estimate:
\begin{eqnarray*}
\hat f_\infty&=&\sum_{j=-1}^{J_\eps} \sum_{\xi\in
\chi_j }\widehat{\beta}_{j, \xi}\ONE_{\{|\widehat{\beta}_{j, \xi
}|\|\psi_{j,\xi}\|_\infty\ge\kappa2^{jd}c_\eps\}}\psi_{j,\xi},
\\
2^{J_\eps d}&=&  \bigl(\eps\sqrt{\log1/\eps} \bigr)^{-1}.
\end{eqnarray*}
Then, for this estimate, we have the following results:

\begin{theorem}
\label{mMinfty}
For $0< r\le\infty$, $\pi\ge1$, $s>\frac{d+1}\pi$,
there exist some constants $c_\infty=c_\infty(s,\pi,r,M)$ such that
if $\kappa^2\ge4 \tau_\infty$, where $\tau_\infty:=\sup_{j, \xi
}2^{-j\fracb{d+1}2}\|\psi_{j,\xi}
\|_\infty$,
\[
\sup_{f\in B^s_{\pi, r}(M)}\bE\|\hat f_\infty-f\|_\infty \le
c_\infty \bigl(\epsilon\sqrt{\log
1/\epsilon} \bigr)^{\frace{s-(d+1)/\pi}{s+d-(d+1)/\pi}}.
\]
\end{theorem}

The following theorem states lower bounds for the minimax rates over
Besov spaces in this model.

\begin{theorem}\label{lowerb} Let $\cEstim$ be the set of all estimators,
for $0< r\le\infty$, $\pi\ge1$, $s>\frac{d+1}\pi$.\vspace*{3pt}
\begin{enumerate}[(b)]
\item[(a)] There exists some constant $C_\infty=C_\infty(s,\pi,r,M)$ such that,
\[
 \inf_{f^\star\in\cEstim}\sup_{f\in B^s_{\pi, r}(M)}\bE\|f^\star
-f\|_\infty\ge
C_\infty \bigl(\epsilon\sqrt{\log
1/\epsilon} \bigr)^{\frace{s-(d+1)/\pi}{s+d-(d+1)/\pi}}.
\]
\item[(b)] For $1\le p< \infty,$ there exists some constant $C_p=C_p(s,\pi
,r,M)$ such that if $s>(\frac{d+1}\pi-\frac{d+1}p)_+$,\vspace*{3pt}
\begin{enumerate}[(3)]
\item[(1)]
If $\frac1p <\frac d{d+1}$
\begin{eqnarray*}
&&\inf_{f^\star\in\cEstim}\sup_{f\in B^s_{\pi, r}(M)} (\bE
\|f^\star-f\|_p^p )^{\fraca1p}\\
&& \quad  \geq C_p
\epsilon^{\frace{s- (d+1)(1/\pi-1/p)}{s +d -(d+1)/\pi}\wedge
\fracc s{s+ d-1/2}\wedge\frace{s-2(1/\pi-1/p)}{s+d-2/\pi}}.
\end{eqnarray*}
\item[(2)]
If $ \frac d{d+1} \leq\frac{1}{p}$ and $d>2$ or $p>1$
\[
\inf_{f^\star\in\cEstim}\sup_{f\in B^s_{\pi, r}(M)} (\bE\|
f^\star-f\|_p^p )^{\fraca1p}
\geq C_p
\epsilon^{\fracc s{s+
d-1/2}\wedge\frace{s-2(1/\pi-1/p)}{s+d-2/\pi}}.
\]
\item[(3)]
If $d=2$ and $p=1$
\[
\inf_{f^\star\in\cEstim}\sup_{f\in B^s_{\pi, r}(M)} (\bE\|
f^\star-f\|_1 )
\geq C_p \epsilon^{\fracc s{s+ 2-1/2}}.
\]
\end{enumerate}
\end{enumerate}
\end{theorem}

\begin{rem}
A careful look at the proof shows that the different rates observed in
the two preceding theorems can be ``explained'' by geometrical
considerations. In fact, depending on the cubature points around which
they are centered, the needlets do not behave the same way. In
particular, their $\bL_p$ norms differ. This leads us to consider two
different regions on the sphere, one near the pole and one closer to
the equator. In these two regions, we considered dense and sparse cases
in the usual way.
This yielded four rates. Then it appeared that one of them (sparse) is
always dominated by the others.\looseness=1
\end{rem}

\section{Applications to the fan beam tomography}\label{sec5}

\subsection{The 2D case: Fan beam tomography}\label{sec5.1}

\begin{figure}
\includegraphics{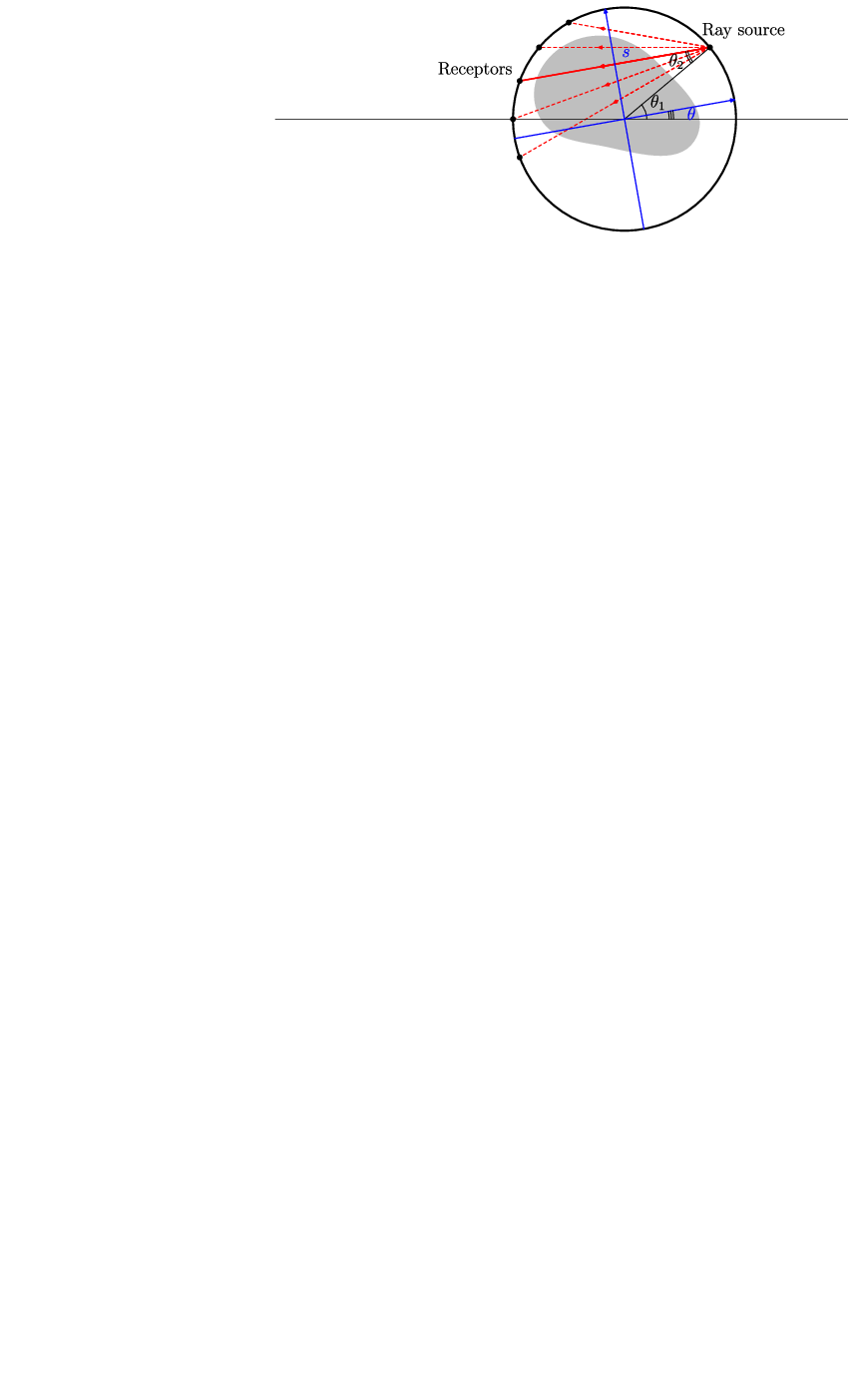}%
\vspace*{-3pt}
\caption{Simplified CAT device.}
\label{fig:CAT}
\vspace*{-3pt}
\end{figure}

When $d=2$, the Radon transform studied in this paper is the fan beam
Radon transform used in a computed axial tomography (CAT) scan.
The geometry of such a device is illustrated in Figure~\ref{fig:CAT}.
An object is placed in the middle of the scanner and X-rays are sent
from a pointwise source, $S(\theta_1)$, making an angle $\theta_1$ with
a reference direction. Rays go through  the object and the energy decay
between the source and an array of receptors is measured.\vadjust{\goodbreak} As the log
decay along the ray is proportional to the integral of the density~$f$
of the object
along the same ray, the measurements are
\[
\tilde Rf(\theta_1,\theta_2) = \int_{e_{\theta_1} + \lambda
e_{(\theta_1-\theta_2)} \in B^2
} f(x)\,\mathrm{d}\lambda
\]
with $e_\theta=(\cos\theta,\sin\theta)$.
This is equivalent to the classical Radon transform
\[
Rf(\theta, s) =
  \int _{\mathop{y\in\theta^\perp,}\limits_{
s\theta+y\in B^2}} f(s\theta + y)\,\mathrm{d}y,
\qquad\theta\in\mathbf{S}^{1},   s \in[-1,1],
\]
for $\theta=\theta_1-\theta_2$ and $s=\sin\theta_2$.
The ray source is then rotated to a different angle and the
measurement process is repeated. In our Gaussian white noise model, we
measure the
continuous function
$Rf(\theta,s)$ through the process $\mathrm{d}Y=Rf(\theta,s)\,\mathrm{d}\theta
\frac{\mathrm{d}s}{(1-s^2)} + \epsilon \,\mathrm{d}W(\theta,s)$. The
measure $\mathrm{d}\theta\frac{\mathrm{d}s}{(1-s^2)}$ corresponds to
the uniform measure $\mathrm{d}\theta_1 \,\mathrm{d}\theta_2$ by the change of
the variable that maps $(\theta_1,\theta_2)$ into $(\theta,s)$.
Our goal is to recover the unknown function, $f,$ from the observation
of $Y$ using the needlet thresholding mechanism described in the
previous sections.

In our implementation, we exploit the tensorial structure of the SVD
basis of the disk in polar coordinates:
\[
f_{k,l,m} (r,\theta)
= (2k+2)^{1/2}P_j^{(0,  l)} (2|r|^2-1)|r|^l Y_{l, m}(\theta), \qquad
   0\leq l \leq k,   k-l =2j,    1\leq m \leq2,
\]
where $P_{j}^{0,l}$ is the corresponding Jacobi polynomial, and
$Y_{l,1}(\theta)=c_l \cos(l\theta)$ and $Y_{l,2}(\theta)=c_l
\sin(l\theta)$
with $c_0=\frac{1}{\sqrt{2\uppi}}$ and $c_l=\frac{1}{\sqrt{\uppi}}$,
otherwise.
The basis of $S^2\times[-1,1]$ has a similar tensorial structure as it
is given by
\[
g_{k,l,m}(\theta, s)
= [h_k]^{-1/2}(1-s^2)^{1/2} C^{1}_k(s) Y_{l,m }(\theta),
\qquad k\ge0,    l \ge0,    1\leq m \leq2,
\]
where $C^{1}_k$ is the Gegenbauer of parameter $1$ and degree $k$.
We recall that the corresponding eigenvalues are
\[
\lambda_k = \frac{2\sqrt{\uppi}}{\sqrt{k+1}}.
\]

\subsection{SVD, needlet and cubature}\label{sec5.2}

In our numerical studies, we compare four different type of
estimators:
linear SVD estimators, thresholded SVD estimators, linear needlet
estimators and thresholded needlet  estimators. They are defined
from the measurement of the values of the Gaussian field on the
SVD basis function $Y_{g_{k,l,m}}$ and the following linear estimates
of, respectively, the SVD basis
coefficients $\langle f,
f_{k,l,m} \rangle$ and the needlet coefficients $\langle f,
\psi_{j,\xi} \rangle$,
\[
\widehat{\alpha}_{k,l,m}   =  \frac{1}{\lambda_k} Y_{
g_{k,l,m}} = \frac{1}{\lambda_k} \int g_{k,l,m}\,\mathrm{d}Y
\]
and
\[
\widehat{\beta}_{j,\xi}  = \sqrt{\lambdaquadra_{j,\xi}} \sum_k
\sqrt{b(k/2^j)}
\sum_{l,m} g_{k,l,m}(\xi) \widehat{\alpha}_{k,l,m}.
\]
The estimators we consider are respectively
defined as:
\begin{eqnarray*}
&&\mbox{linear SVD estimates}\hspace*{55pt}
\hat{f}^{\mathrm{LS}}_J  = \sum_{k<2^J} \sum_{l,m} \widehat{\alpha}_{k,l,m}
f_{k,l,m},\\
&&\mbox{linear needlet estimates}\hspace*{44pt}
\hat{f}^{\mathrm{LN}}_J  = \sum_{j<J} \sum_{\xi} \widehat{\beta}_{j,\xi}
\psi_{j,\xi},\\
&&\mbox{thresholded SVD estimates}\hspace*{31pt}
\hat{f}^{\mathrm{TS}}_{T}  = \sum_{k<2^J} \sum_{l,m} \rho_{T_k} (
\widehat{\alpha}_{k,l,m}
 ) f_{k,l,m},\\
&&\mbox{thresholded needlet estimates} \qquad
\hat{f}^{\mathrm{TN}}_{T}  = \sum_{j<J} \sum_{\xi}
\rho_{T_{j,\xi}} (\widehat{\beta}_{j,\xi} ) \psi
_{j,\xi},
\end{eqnarray*}
where $\rho_T(\cdot)$ is the hard threshold function with threshold
$T$:
\[
\rho_T(x) =
\cases{\displaystyle
x ,&\quad  if $|x|\geq T$,\cr\displaystyle
0 ,&\quad  otherwise.
}
\]
A more precise description is given in Table~\ref{tab:algo}.
In our experiments, the values of $Y_{g_{k,l,m}}$ have been obtained
from an initial approximation of $\langle f, f_{k,l,m}\rangle$
computed with a very fine cubature to which a Gaussian i.i.d. sequence
is added.

\begin{table*}
\tabcolsep=0pt
\caption{Algorithmic description of the considered estimators}
\label{tab:algo}
 \begin{tabular*}{\textwidth}{@{\extracolsep{4in minus 4in}}lcccc@{}}
 \hline
& Linear SVD & Thresh. SVD & Linear needlet & Thresh. needlet \\
& $\hat{f}^{\mathrm{LS}}$ & $\hat{f}^{\mathrm{TS}}$ & $\hat{f}^{\mathrm{LN}}$ & $\hat
{f}^{\mathrm{TN}}$\\
\hline
Observation& \multicolumn{4}{c@{}}{$  \mathrm{d}Y=Rf \,\mathrm{d}\theta\frac{\mathrm{d}s}{\sqrt{1-s^2}} + \epsilon
\,\mathrm{d}W$}\\
& \multicolumn{4}{c@{}}{\hrulefill}\\
SVD Dec&\multicolumn{4}{c@{}}{$ Y_{h_{k,l,}}= \langle Y, g_{k,l,m} \rangle= \langle Rf ,
g_{k,l,m} \rangle+ \epsilon_{k,l,m} = \mu_k \langle f ,
f_{k,l,m} \rangle+ \epsilon_{k,l,m}$}\\
& \multicolumn{4}{c@{}}{\hrulefill}\\
Inv. Radon& \multicolumn{4}{c@{}}{$\widehat{\alpha}_{k,l,m}=\frac{1}{\mu_k} \langle Y , g_{k,l,m} \rangle=
\langle f , f_{k,l,m}\rangle+ \frac{1}{\mu_k}\epsilon_{k,l,m}$} \\
& \multicolumn{4}{c@{}}{\hrulefill}\\
 \multirow{2}{31pt}{Needlet \mbox{\quad transf.}}&&& \multicolumn{2}{c@{}}{\hspace*{-21pt}$\widehat{\beta}_{j,\xi}=
\sqrt{\lambdaquadra_{j,\xi}} \sum_{k}\sqrt{b (\frac{k}{2^j} )}$}\\
&&&\multicolumn{2}{c@{}}{\hspace*{25pt}${}\times\sum_{l,m} f_{k,l,m}(\xi) \widehat{\alpha}_{k,l,m}$}
\\
&&& \multicolumn{2}{c@{}}{\hrulefill}\\
Coeff. mod.& $\widehat{\alpha}_{k,l,m}^{SL}=\mathbf{1}_ {k\leq k_{\max}}\widehat{\alpha}_{k,l,m}$
& $\widehat{\alpha}_{k,l,m}^{ST}=\rho_{T_k} (\widehat{\alpha}_{k,l,m} )$
& $\widehat{\beta}_{j,\xi}^{NL}= \mathbf{1}_{j<j_{\max}}\widehat{\beta}_{j,\xi}$
& $\widehat{\beta}_{j,\xi}^{NT}= \rho_{T_{j,\xi}} (\widehat{\beta}_{j,\xi} )$
\\
&\multicolumn{1}{c}{\hrulefill}&\multicolumn{1}{c}{\hrulefill}&\multicolumn{1}{c}{\hrulefill}&\multicolumn{1}{c@{}}{\hrulefill}\\[-1pt]
\multirow{2}{43pt}{Needlet inv.}& &
& \multicolumn{2}{c@{}}{\hspace*{-66pt}$\widehat{\alpha}_{k,l,m}^{\star}=\sum_{j} \sqrt{b (\frac{k}{2^j}
)}$}\\
& &
& \multicolumn{2}{c@{}}{\hspace*{35pt}${}\times\sum_{\xi\in\chi_j} \sqrt{\lambdaquadra_{j,\xi}}f_{k,l,m}(\xi) \widehat{\beta}_{j,\xi}^{\star}$}
\\
&&& \multicolumn{2}{c@{}}{\hrulefill}\\
SVD rec.&\multicolumn{4}{c@{}}{$\hat f^{\star}=
\sum_{k,l,m} \widehat{\alpha}_{k,l,m}^{\star}f_{k,l,m}$}
\\
\hline
\end{tabular*}
\end{table*}

We have used, in our numerical experiments, thresholds of the form
\[
T_k = \frac{\kappa}{\lambda_k}  \epsilon\sqrt{\log1/
\epsilon} \quad\mbox{and}\quad
T_{j,\xi} = \kappa \sigma_{j,\xi}  \epsilon\sqrt{\log1/
\epsilon},
\]
where $\sigma_{j,\xi}$ is the standard deviation of the
noisy needlet coefficients when $f=0$ and~\mbox{$\epsilon=1$}:
\[
\sigma_{j,\xi}^2= \lambdaquadra_{j,\xi} \sum_k b(k/2^j)
\sum_{l,m} g_{k,l,m}(\xi)^2.
\]
Note that, while the needlet threshold is different than in
Theorem~\ref{mM}, as $\sigma_{j,\xi}$ is of order~$2^{j\nu}$,
its conclusions remain valid.

An important issue in the needlet scheme is the choice of the
cubature in the needlet construction. Proposition~\ref{prop:CUB}\vadjust{\goodbreak}
ensures the existence of a suitable cubature $\xi_j$ for every
level~$j$ based on a cubature $\tilde\xi_j$ on the sphere
but gives neither an explicit construction of the points
on the sphere nor an explicit formula for the weights $\lambdaquadra
_{j,\xi}$.
Those ingredients are, nevertheless, central in the numerical scheme and
should be specified.

Three possibilities have been considered:
a numerical cubature deduced from an almost uniform cubature of the
half sphere available, an approximate cubature deduced from the
Healpix cubature on the sphere and a cubature obtained by subsampling
a tensorial cubature associated to the latitude and longitude
coordinates on the sphere. The first strategy has been considered by
Baldi  \textit{et al.}  \cite{baldi-2008} in a slightly different context;
there is, however, a strong limitation on the maximum
degree of the cubature available and, thus, this solution has been
abandoned. The Healpix strategy, also considered by Baldi  \textit{et al.}  in
another paper \cite{baldi09asymp},
is easily implementable but, as it is based on an approximate
cubature, fails to be precise enough in our numerical experiments.
The last strategy relies on the
subsampling on a tensorial grid on the sphere. While such a strategy
provides a simple way to construct an admissible cubature, the
computation of the cubature weights is becoming an issue as no closed
form is available.

To overcome those issues, we have considered a cubature formula based
on the full tensorial grid appearing as proposed by \cite
{muciaccia97fastspherharmonanaly}. This
cubature does not satisfy the condition of
Proposition~\ref{prop:CUB}, but its weights can be computed explicitly.
Furthermore,
we argue that, using our modified threshold, we can still control the
risk of the estimator. Indeed, note first that the modified threshold
is such
that the thresholding of a needlet depends only on its scale parameter
$j$ and on
its center $\xi$ and not on the corresponding cubature weight~$\lambdaquadra_{j,\xi}$.
Assume now that we have a collection of $K$ cubature, each
satisfying conditions of Proposition~\ref{prop:CUB} and thus defining
a suitable estimate, $\hat{f}_k$. We can use the ``average''
cubature obtained by adding all the cubature points and averaging the
cubature weights. This new
cubature defines a new estimate, $\hat{f}$, satisfying
\[
\hat{f} = \frac{1}{K} \sum_{k=1}^K \hat{f}_k.
\]
By convexity, for any $p\geq1$,
\[
\|f-\hat{f}\|_p^p = \Biggl\| f - \frac{1}{K} \sum_{k=1}^K \hat{f}_k\Biggr\|_p^p
  \leq \frac{1}{K} \sum_{k=1}^K
\| f - \hat{f}_k\|_p^p
\]
and, thus, this average estimator is as efficient as the worst estimator
in the family~$\hat{f}_k$. The full tensorial cubature
is an average of suitable cubatures. The corresponding
estimator satisfies the error bounds of
Theorems~\ref{mM} and~\ref{mMinfty}.
Note that this principle is quite close to the cycle-spinning method
introduced by Donoho  \textit{et al.}  Indeed, the same kind of numerical gain is
obtained with this method.
The numerical comparison of the Healpix cubature and our tensorial
cubature is largely in favor of our scheme. Furthermore, as already
noticed by \cite{muciaccia97fastspherharmonanaly}, the tensorial
structure of
the cubature leads to some simplification in the numerical
implementation of the needlet estimator. The resulting scheme is almost
as fast as the Healpix-based one.

\begin{figure}

\includegraphics{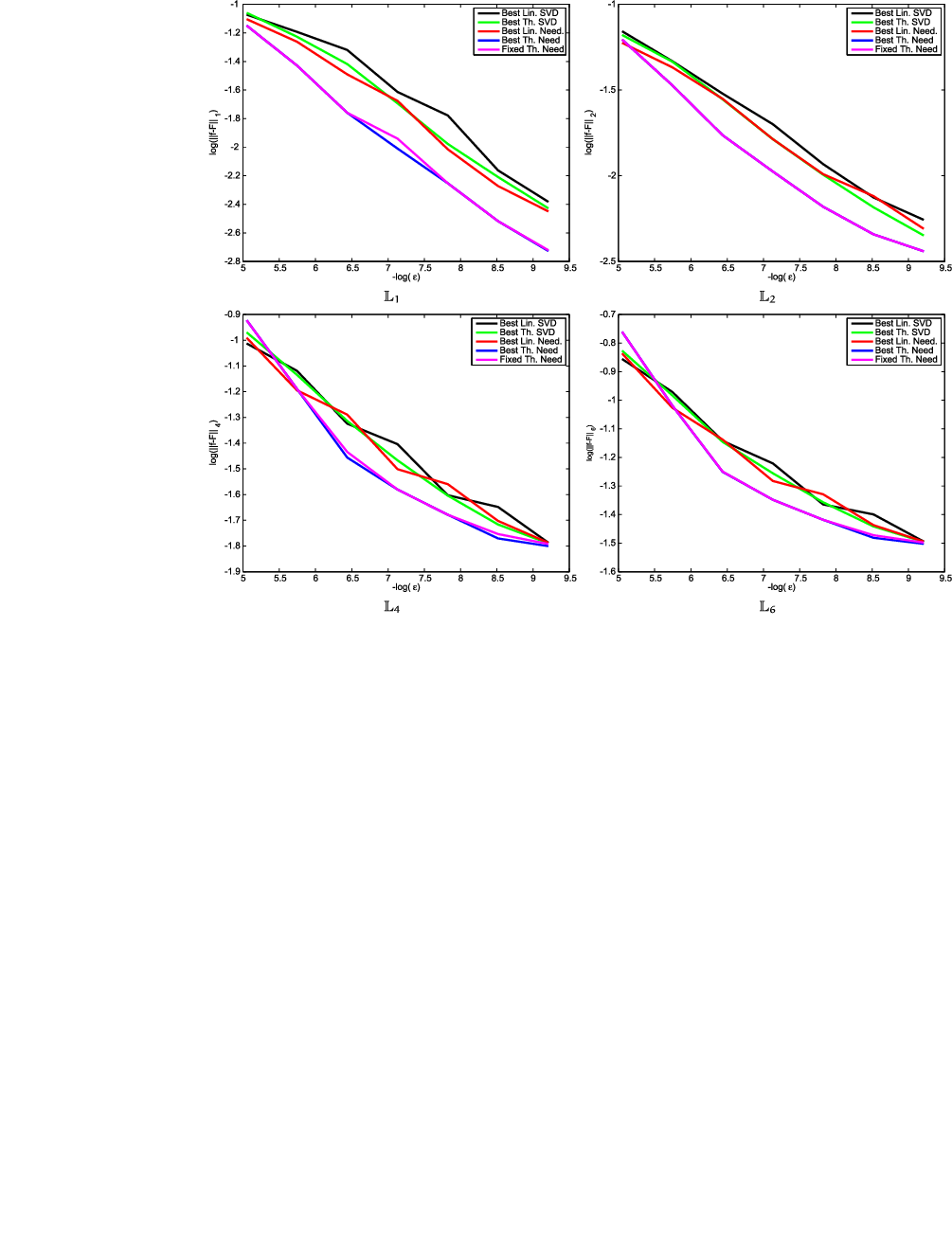}

\caption{Estimation results in log $\bL_p$ norm. Each figure shows the
decay of the logarithm of the error against the logarithm of the
noise parameter for the specified norm.}
\label{fig:result}
\end{figure}

\subsection{Numerical results}\label{sec5.3}

\begin{figure}

\includegraphics{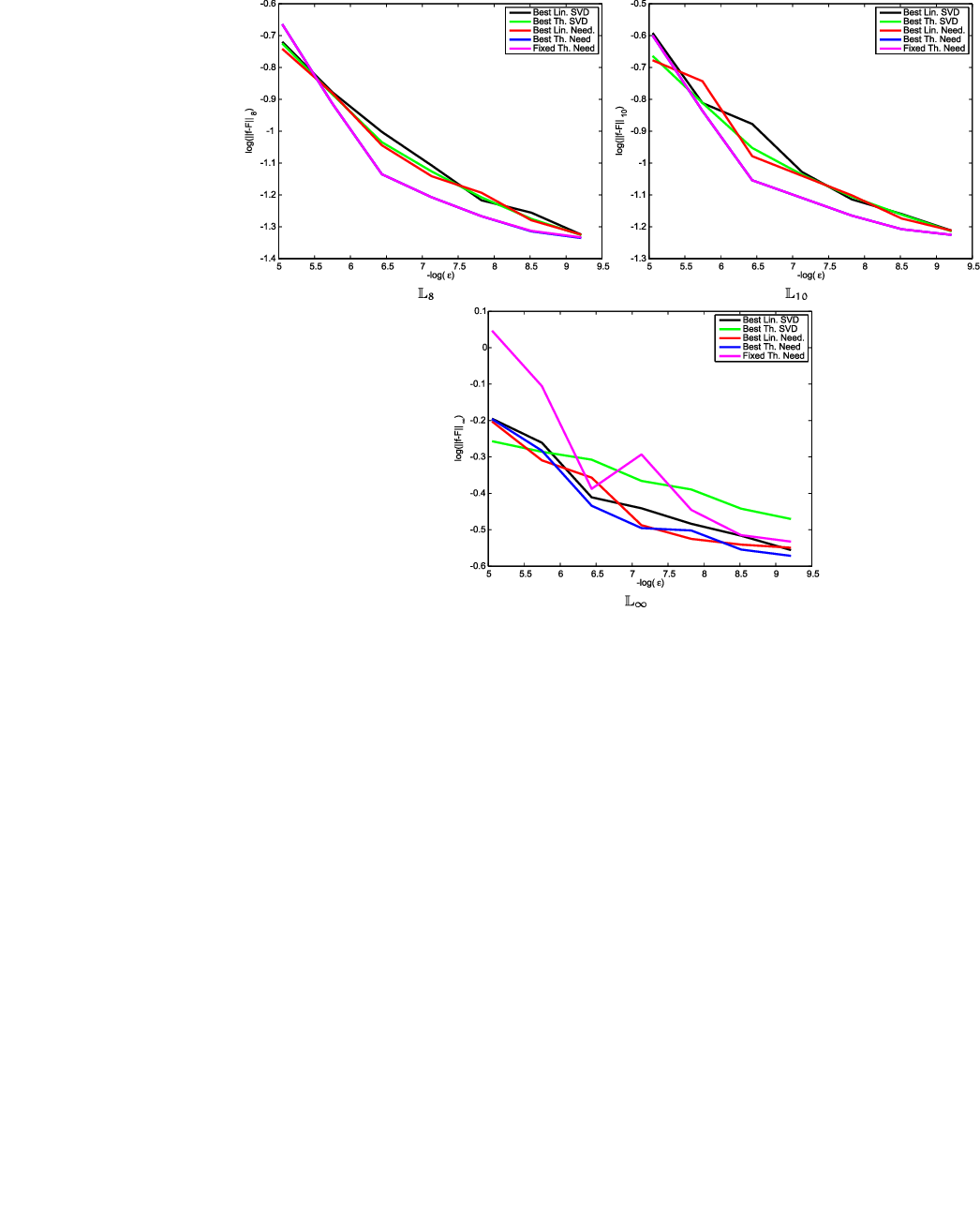}

\caption{Estimation results in log $\bL_p$ norm. Each figure shows the
decay of the logarithm of the error against the logarithm of the
noise parameter for the specified norm.}
\label{fig:result2}
\end{figure}

\begin{figure}

\includegraphics{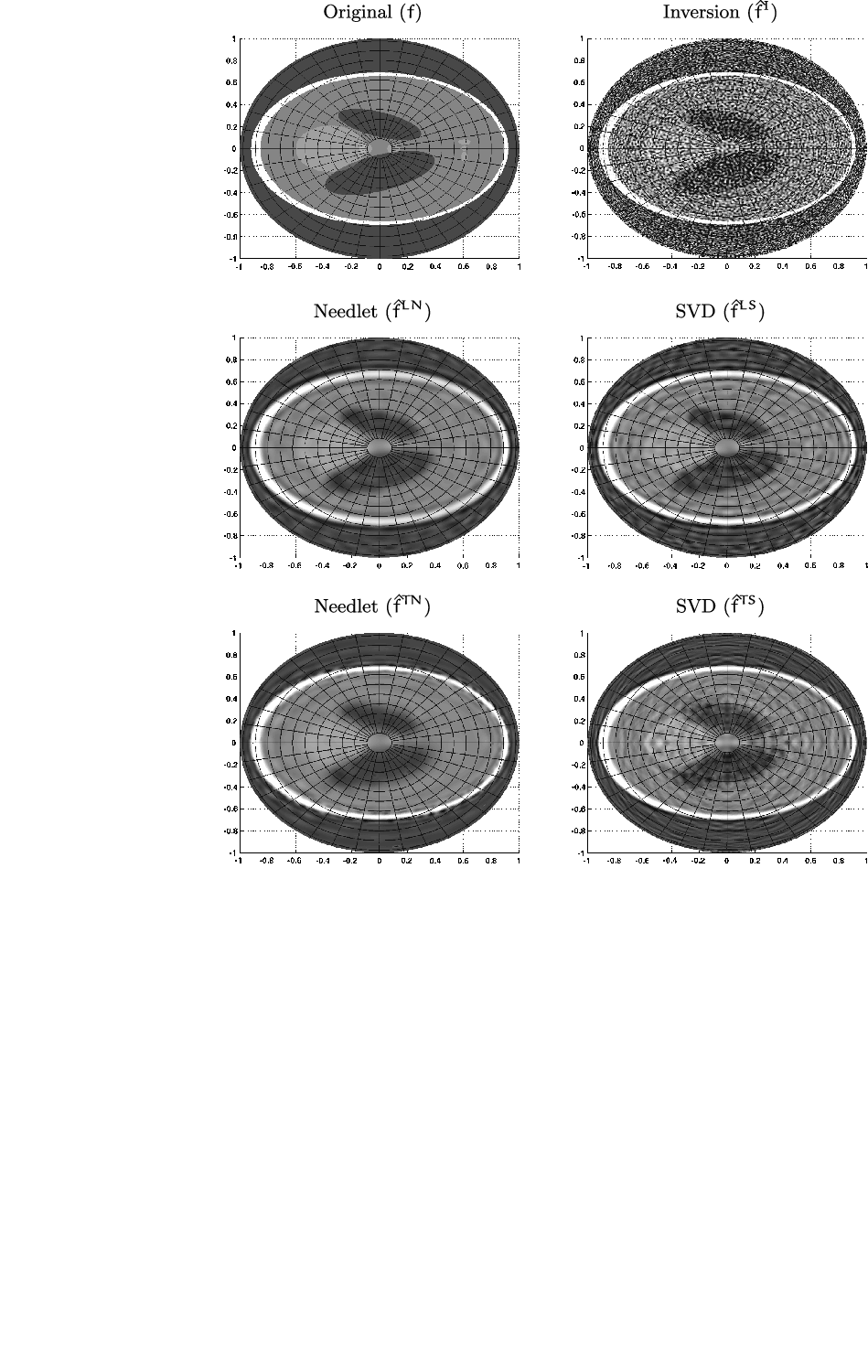}

\caption{Visual comparison for the original Logan--Shepp phantom with
$\epsilon=8/1\,000$. Errors are much more localized in the
needlet-based estimates compared to the fully delocalized errors of
the SVD based estimates. Fine structures are much more restored in
the thresholded needlet estimate than in the other estimates.}
\label{fig:numcomp}
\end{figure}

In this section, we compare five ``estimators'' (linear SVD with best
scale, linear needlet with best scale,
thresholded SVD with best $\kappa$, thresholded needlet with best
$\kappa$ and
thresholded needlet with $\kappa=3$)
for seven different norms ($\bL_1$, $\bL_2$, $\bL_4$, $\bL_6$, $\bL
_7$, $\bL_8$,
$\bL_{10}$ and $\bL_\infty$) and seven
noise levels $\epsilon$ ($2^k/1\,000$ for $k$ in $0,1,\ldots,6$).
Each subfigure of Figures~\ref{fig:result} and~\ref{fig:result2}
plots the logarithm of the estimation error
for a specific norm against the opposite of the logarithm of the
noise level. Note that the subfigure overall aspect is
explained by the error decay when the noise level diminishes. The
good theoretical behavior of the thresholded needlet estimator is
confirmed numerically:
the thresholded needlet estimator with an optimized $\kappa$ appears
as the best estimator for every norm while a fixed $\kappa$
yields a~very good estimator, except for the $\bL_\infty$ case, as
expected by our theoretical results. These results are confirmed
visually by the
reconstructions of Figure~\ref{fig:numcomp}. In the needlet ones,
errors are smaller and much more localized than in their
SVD counterparts. Observe also how the fine structures are much more
preserved with the thresholded needlet estimate than with any other methods.

We conclude this paper with some sections devoted to the proofs of our results.

\section{Needlet properties}\label{besov}\label{sec6}
\subsection{Key inequalities}\label{sec6.1}
The following inequalities are true (and proved in \cite
{pnarco,NPW,PXU,pxukball}) and will be fundamental in the sequel.
In the following lines, $g_{j,\xi}$ will stand either for $\phi
_{j,\xi}$ or $\psi_{j,\xi}$:
%
\begin{eqnarray}
\label{A}
&\displaystyle\forall j\in\bN,\ \forall \xi\in\chi_j \qquad0< c\leq\|
g_{j,\xi}\|_2^2 \leq1,&\vadjust{\goodbreak}
\\
&\displaystyle\forall j\in\bN, \xi\in\chi_j,  \ \forall x \in\cX  \qquad
\sum_{\xi\in\chi_j}  \| g_{j,\xi}\|_1| g_{j,\xi}(x) | \leq C
<\infty,&
\\
\label{EqFond2}
&\displaystyle|g_{j,\xi}(x)|\le C_M\frac{2^{jd/2}}{\sqrt{W_j(x)}(1+2^jd(x,\xi))^M}&
\end{eqnarray}
(recall that $W_j(x)$ has been defined in (\ref{Wj})). From these
inequalities, one can deduce the following ones (see \cite{kyoto}):
for all $1 \leq p \leq\infty$,
%
\begin{eqnarray}
\label{1}
 \biggl( \sum_{\xi\in\chi_j} |\langle f, g_{j,\xi}\rangle|^p \|
g_{j,\xi} \|^p_p \biggr)^{1/p}
&\leq& C \| f\|_p,
\\[-2pt]
\label{2}
 \biggl\| \sum_{\xi\in\chi_j} \lambda_\xi g_{j,\xi}(x) \biggr\|_p &\leq&
\biggl(\frac Cc \biggr)^2 \biggl( \sum_{\xi\in\chi_j} \|\lambda_\xi
g_{j,\xi} \|^p_p \biggr)^{1/p}.
\end{eqnarray}

\subsection{Besov embeddings}\label{sec6.2}
It is a key point to clarify how the Besov spaces defined above may be
included in each of the others. As will be seen, the embeddings\vadjust{\goodbreak} will
parallel the standard embeddings of usual Besov spaces, but with
important differences, which, in particular, yield new minimax rates of
convergence as detailed above.

We begin with an evaluation of the different $\bL_p$ norms of the needlets.
More precisely, in \cite{pxukball} it is shown that, for $0<p \le
\infty,$
%
\begin{equation}
\label{norm-est}
\|\psi_{j,\xi}\|_p \sim\|\phi_{j,\xi}\|_p
\sim \biggl(\frac{2^{jd}}{W_j(\xi)}  \biggr)^{1/2-1/p},
\qquad\xi\in\chi_j.
\end{equation}
The following inequalities are proved in \cite{5authors}:
%
\begin{eqnarray}
\label{B}
\sum_{\xi\in\chi_j} \| g_{j,\xi} \|^p_p &\leq& c2^{j(dp/2 + (p/2 - 2)_+)}
\qquad\mbox{if }  p\neq4,
\\
\label{BI}
\sum_{\xi\in\chi_j} \|g_{j,\xi} \|^p_p &\leq& cj 2^{jdp/2}
\qquad\mbox{if }  p=4.
\end{eqnarray}
 We are now able to state the embedding results (see \cite{kyoto}).

\begin{theorem} \label{embb}
(1)
$ { 1\leq p\leq\pi \leq\infty \Rightarrow
B^{s}_{\pi,r} \subseteq B^{s}_{p,r}}.$

(2)
\(
\infty\geq p\geq\pi >0 $,   $s >(d+1) (1/\pi- 1/p)$,
$\Rightarrow
B^{s}_{\pi,r} \subseteq B^{s-(d+1) (1/\pi- 1/p)}_{p,r}.
\)
\end{theorem}

\section{Proof of the upper bounds}\label{sec7}
A important tool for the proof of the upper bounds that clarifies the
thresholding procedure is the following lemma.
\begin{lemma}\label{variance}
For all $j\ge-1$, $\xi\in\chi_j$,
$\widehat\beta_{j,\xi}$ has a Gaussian distribution with mean $\beta
_{j,\xi}$ and variance $\sigma^2_{j,\xi}$, with
\[
\sigma^2_{j,\xi}\le c2^{j(d-1)} \eps^2.
\]
\end{lemma}

\begin{pf}
As we can write
\begin{eqnarray*}
\widehat{\beta}_{j, \xi}
&=&\sum_{k, l, m} \gamma^{j,\xi}_{k,l,m}\int_{B^d}ff_{k,l,m}\,\mathrm{d}x
+ \sum_{k, l, m}\gamma^{j,\xi}_{k,l,m}\frac{\eps}{\lambda_k}
Z_{k,l,m}\\
&=&\beta_{j, \xi} + Z_{j,\xi}.
\end{eqnarray*}
Here the summation is over
$\{(k,l,m)\dvt 0\le k < 2^j, 0\le l \le k, l\equiv k\ (\operatorname{mod} 2) ,1 \leq m
\leq N_{d-1}(l)\}$.
Since the $Z_{k,l,m} $'s are independent $N(0,1)$ random variables,
$ Z_{j,\xi} \sim N(0,\sigma^2_{j,\xi}),$
we have
%
\begin{equation}\label{sigma-j}
\sigma^2_{j,\xi}
= \eps^2 \sum_{k, l, m}
|\gamma^{j,\xi}_{k,l,m}|^2 \frac{(k)_d}{\uppi^{d-1}2^d k}
\le\frac{(2^j)_{d-1}}{ \uppi^{d-1} 2^d}
\leq c2^{j(d-1)} \eps^2
\end{equation}
with $c=(d/2\uppi)^{d-1}$.
Here, we used that $\{f_{k,l,m}\}$ is an orthonormal basis for $\bL_2$
and hence
$\sum_{k, l, m} |\gamma^{j,\xi}_{k,l,m}|^2 = \|\psi_{j, \xi}\|
_2^2\le1$.\
\end{pf}

Let us now begin with the second theorem, the proof of which is
slightly simpler.

\subsection{\texorpdfstring{Proof of Theorem \protect\ref{mMinfty}}{Proof of Theorem 2}}
\label{sec7.1}

In this proof, as well as in the other one, $C$ will denote any
\textit{constant} in the sense of Theorems~\mbox{\ref{mM}--\ref{lowerb}}; we assume that $f$ belongs to the Besov ball
$B^s_{\pi,r}(M)$ and
measure the loss in $\bL_p$ norm (here, $p=\infty$). Then $C$ denotes
a generic
quantity that depends only on $s$, $r$, $p$ and $M$. Note that the
exact value denoted by $C$ may vary from one line to another.

We have, if we denote
\begin{eqnarray*}
A_J (f)&:=& \sum_{j>J} \sum_{\xi\in
\chi_j }{\beta}_{j, \xi}\psi_{j,\xi},
\\
\|\hat{f}_\infty-f\|_\infty
&\leq&
\|\hat{f}_\infty-
A_J(f)\|_\infty+\|A_J(f) -f\|_\infty\\
&\leq&
\|\hat{f}_\infty- A_J(f)\|_\infty + C \| f \|_{B^s_{\pi,r} }
2^{-J(s-(d+1)/\pi)}.
\end{eqnarray*}
We have used, as $B^s_{\pi,r} \subset B^{s- (d+1)\fraca1\pi}_{\infty,r}$,
\[
\| A_Jf -f \|_\infty\leq C \| f \|_{B^s_{\pi,r} }
2^{-J(s- (d+1)\fraca1\pi)}.
\]
Moreover,
\[
2^{-J(s-(d+1)/\pi)} \leq  \bigl(\epsilon\sqrt{\log1/\epsilon
} \bigr)^{-\fracb{s-(d+1)/\pi}{d}}
\leq  \bigl(\epsilon\sqrt{\log1/\epsilon} \bigr)^{-\frace
{s-(d+1)/\pi}{ s-(d+1)/\pi+d}}
\]
as $s>\frac{d+1}\pi$.

We have, using \eqref{2},
\begin{eqnarray*}
\|\hat{f}_\infty- A_J(f)\|_\infty
&\leq&
\sum_{j<J} \biggl\| \sum_{\xi\in
\chi_j}  \bigl( \widehat{\beta}_{j,\xi}
\Charac{|\widehat{\beta}_{j,\xi} |\|\psi_{j,\xi}\|_\infty\geq\kappa
2^{jd} \epsilon\sqrt{\log1/\epsilon}}
-{\beta}_{j,\xi}  \bigr)\psi_{j,\xi} \biggr\|_\infty
\\
&\le&\sum_{j<J} \biggl( \biggl\| \sum_{\xi\in\chi_j}  \bigl( (\widehat{\beta
}_{j,\xi} -\beta_{j,\xi} )\psi_{j,\xi}
\Charac{|\widehat{\beta}_{j,\xi} |\|\psi_{j,\xi}\|_\infty\geq\kappa
2^{jd} \epsilon\sqrt{\log1/\epsilon} }
 \bigr) \biggr\|_\infty
\\
&&\hphantom{\sum_{j<J} \biggl(}
{}   + \biggl\| \sum_{\xi\in\chi_j}  \bigl( \beta_{j,\xi}
\psi_{j,\xi}
\Charac{|\widehat{\beta}_{j,\xi} |\|\psi_{j,\xi}\|_\infty< \kappa
2^{jd} \epsilon\sqrt{\log1/\epsilon}}
 \bigr) \biggr\|_\infty \biggr)\\
&\leq&\frac Cc
\sum_{j<J}  \biggl( \sup_{\xi\in\chi_j}  \bigl(|\widehat{\beta
}_{j,\xi} -\beta_{j,\xi} | \|\psi_{j,\xi}\|_\infty
\Charac{|\widehat{\beta}_{j,\xi} |\|\psi_{j,\xi}\|_\infty\geq\kappa
2^{jd} \epsilon\sqrt{\log1/\epsilon} }  \bigr)
\\ &&\hphantom{\frac Cc
\sum_{j<J}  \biggl(}{}    + \sup_{\xi\in\chi_j} \bigl ( |\beta_{j,\xi}|
\| \psi_{j,\xi}\|_\infty
\Charac{|\widehat{\beta}_{j,\xi} |\|\psi_{j,\xi}\|_\infty< \kappa
2^{jd} \epsilon\sqrt{\log1/\epsilon} } \bigr)
 \biggr).
\end{eqnarray*}
We decompose the first term of the last inequality
\begin{eqnarray*}
&& \sup_{\xi\in\chi_j}  \bigl(|\widehat{\beta}_{j,\xi} -\beta_{j,\xi
} | \|\psi_{j,\xi}\|_\infty
\Charac{|\widehat{\beta}_{j,\xi} |\|\psi_{j,\xi}\|_\infty\geq\kappa
2^{jd} \epsilon\sqrt{\log1/\epsilon}} \bigr)
\\[-2pt]
&& \quad = \sup_{\xi\in\chi_j}  \bigl(|\widehat{\beta}_{j,\xi} -\beta
_{j,\xi} | \|\psi_{j,\xi}\|_\infty
\Charac{|\widehat{\beta}_{j,\xi} |\|\psi_{j,\xi}\|_\infty\geq\kappa2^{jd}
\epsilon\sqrt{\log1/\epsilon}}\\[-2pt]
&& \hphantom{\sup_{\xi\in\chi_j}  \bigl(}\qquad {}
\times \bigl( \Charac{|\beta_{j,\xi} |\|\psi_{j,\xi}\|_\infty
\geq\fracd{\kappa}2 2^{jd} \epsilon\sqrt{\log1/\epsilon} }
+ \Charac{|\beta_{j,\xi} |\|\psi_{j,\xi}\|_\infty< \fracd{\kappa
}2 2^{jd} \epsilon\sqrt{\log1/\epsilon} }  \bigr) \bigr)
\\[-2pt]
&& \quad  \le
\sup_{\xi\in\chi_j}  \bigl( |\widehat{\beta}_{j,\xi} -\beta_{j,\xi
} | \|\psi_{j,\xi}\|_\infty
\Charac{|\widehat{\beta}_{j,\xi}-\beta_{j,\xi} |\|\psi_{j,\xi}\|
_\infty\geq\fracd{\kappa}2 2^{jd} \epsilon\sqrt{\log1/\epsilon}
} \bigr)
\\[-2pt]
&& \qquad {}  +
\sup_{\xi\in\chi_j}  \bigl(|\widehat{\beta}_{j,\xi} -\beta_{j,\xi}
| \|\psi_{j,\xi}\|_\infty
\Charac{| \beta_{j,\xi} |\|\psi_{j,\xi}\|_\infty>\fracd{\kappa}2
2^{jd} \epsilon\sqrt{\log1/\epsilon}} \bigr)
\end{eqnarray*}
and the second one
\begin{eqnarray*}
&&\sup_{\xi\in\chi_j}  \bigl( |\beta_{j,\xi}| \| \psi_{j,\xi}\|
_\infty
\Charac{|\widehat{\beta}_{j,\xi} |\|\psi_{j,\xi}\|_\infty< \kappa
2^{jd} \epsilon\sqrt{\log1/\epsilon} } \bigr) \\[-2pt]
&& \quad =
\sup_{\xi\in\chi_j}  \bigl( |\beta_{j,\xi}| \| \psi_{j,\xi}\|
_\infty
\Charac{|\widehat{\beta}_{j,\xi} |\|\psi_{j,\xi}\|_\infty< \kappa
2^{jd} \epsilon\sqrt{\log1/\epsilon} } \\[-2pt]
&&\hphantom{\sup_{\xi\in\chi_j}  \bigl(} \qquad {}    \times
 \bigl( \Charac{ | \beta_{j,\xi} | \|\psi_{j,\xi}\|_\infty \geq2
\kappa2^{jd} \epsilon\sqrt{\log1/\epsilon} }
+\Charac{ | \beta_{j,\xi} | \|\psi_{j,\xi}\|_\infty< \kappa
2^{jd} \epsilon\sqrt{\log1/\epsilon} }
 \bigr)  \bigr)\\[-2pt]
 && \quad \leq
\sup_{\xi\in\chi_j}  \bigl( |\beta_{j,\xi}| \| \psi_{j,\xi}\|
_\infty
\Charac{|\widehat{\beta}_{j,\xi} -\beta_{j,\xi} |\|\psi_{j,\xi}\|
_\infty> \kappa2^{jd} \epsilon\sqrt{\log1/\epsilon} } \bigr)
\\[-2pt]
&& \qquad {}  + \sup_{\xi\in\chi_j}  \bigl(|\beta_{j,\xi}| \| \psi
_{j,\xi}\|_\infty
\Charac{| \beta_{j,\xi} |\|\psi_{j,\xi}\|_\infty< \kappa2^{jd}
\epsilon\sqrt{\log1/\epsilon} } \bigr) .
\end{eqnarray*}
Now we will bound each of the four terms coming from the last two inequalities.
Since for $X\sim N(0,\sigma^2),$ we have
\[
\bE\bigl( |Y| \Charac{|Y| > \lambda\sigma} \bigr)
= \sigma \frac2{\sqrt{2\uppi}} \int_\lambda^\infty y \mathrm{e}^{-y^2/2}\,\mathrm{d}y
=\mathrm{e}^{- \lambda^2/2} \frac2{\sqrt{2\uppi}}\le \mathrm{e}^{- \lambda^2/2}
.
\]
Noticing that the standard deviation of $(\widehat{\beta}_{j,\xi} -\beta
_{j,\xi} )
\|\psi_{j,\xi}\|_\infty$ is smaller than
$\tau_\infty2^{jd} \epsilon $ (using Lemma \ref{variance} and
\eqref{EqFond2}),
we have
\begin{eqnarray*}
&&\sum_{j\leq J}\bE \biggl( \sup_{\xi\in\chi_j}  \bigl(
|\widehat{\beta}_{j,\xi} -\beta_{j,\xi} | \|\psi_{j,\xi}\|_\infty
\Charac{|\widehat{\beta}_{j,\xi}-\beta_{j,\xi} |\|\psi_{j,\xi}\|
_\infty\geq\fracd{\kappa}2 2^{jd} \epsilon\sqrt{\log1/\epsilon}
}  \bigr)  \biggr)
\\[-2pt]
&& \quad \leq\sum_{j\leq J} \sum_{\xi\in\chi_j} \bE \bigl(|\widehat{\beta
}_{j,\xi} -\beta_{j,\xi} | \|\psi_{j,\xi}\|_\infty
\Charac{|\widehat{\beta}_{j,\xi}-\beta_{j,\xi} |\|\psi_{j,\xi}\|
_\infty\geq\fracd{\kappa}{2} 2^{jd} \epsilon\sqrt{\log1/\epsilon
} }  \bigr)\leq c2^{Jd} \epsilon^{\kappa^2/2\tau_\infty^2}
\\[-2pt] && \quad
\le C \epsilon^{\kappa^2/2\tau_\infty^2 -1} \sqrt{\log1/\epsilon
}^{-1} \leq C \epsilon\sqrt{\log1/\epsilon}
\end{eqnarray*}
if $\kappa^2 \geq4\tau_\infty^2,$ where we have used $\Card
{\chi_j}\le c2^{jd}$.
This proves that this term will be of the right order.

For the second term,
let us observe that we have, using Theorem \ref{embb},
\[
| \beta_{j,\xi} |\|\psi_{j,\xi}\|_\infty\leq C \| f
\|_{B^s_{\pi,r} } 2^{-j(s-(d+1)/\pi)},\vadjust{\goodbreak}
\]
so only the $j$'s index such that $j\leq j_1$ will verify this inequality:
\[
2^{j_1} \sim  \biggl(\frac{2C \| f \|_{B^s_{\pi,r}
}}{\kappa}  \bigl(\epsilon\sqrt{\log1/\epsilon} \bigr)
\biggr)^{-\fracc
1{s+d-(d+1)\pi}}.
\]
On the other side, using the Pisier lemma \cite{Pisierentropie},
\[
\bE \biggl(\sup_{\xi\in\chi_j}  (|\widehat{\beta}_{j,\xi}
-\beta_{j,\xi} | \|\psi_{j,\xi}\|_\infty ) \biggr)
\leq\tau_\infty2^{jd} \epsilon\sqrt{2 \log{2 c2^{jd}}}.
\]
So
\begin{eqnarray*}
&&\sum_{j\leq j_1 }\bE \Bigl(
\sup_{\xi\in\chi_j} \bigl (|\widehat{\beta}_{j,\xi} -\beta_{j,\xi}
| \|\psi_{j,\xi}\|_\infty
\Charac{| \beta_{j,\xi} |\|\psi_{j,\xi}\|_\infty> \fracd{\kappa
}2 2^{jd} \epsilon\sqrt{\log1/\epsilon} }
 \bigr)  \Bigr) \\
&& \quad \leq\tau_\infty
\sum_{j\leq j_1 }
2^{jd} \epsilon\sqrt{2 \log{2c 2^{jd}}}
\leq C \epsilon j_1^{\fraca12} 2^{j_1d} \\
&& \quad \lesssim C  (\| f \|
_{B^s_{\pi,r} } )^{\frace{s-(d+1)\pi}{s+d-(d+1)/\pi}}
\bigl(\epsilon\sqrt{\log1/\epsilon} \bigr)^{\frace{s-(d+1)/\pi
}{s+d-(d+1)/\pi}}.
\end{eqnarray*}
This proves that this term will be of the right order.
Concerning the first term of the second inequality,
\begin{eqnarray*}&&\bE \Bigl(\sup_{\xi\in\chi_j}  \bigl( |\beta
_{j,\xi}| \| \psi_{j,\xi}\|_\infty
\Charac{|\widehat{\beta}_{j,\xi} -\beta_{j,\xi} |\|\psi_{j,\xi}\|
_\infty> \kappa2^{jd} \epsilon\sqrt{\log1/\epsilon} }
 \bigr)  \Bigr)\\
&& \quad \leq C \| f \|_{B^s_{\pi,r} } \sum_{\xi\in\chi_j} P \bigl( |\widehat
{\beta}_{j,\xi} -\beta_{j,\xi} |\|\psi_{j,\xi}\|_\infty> \kappa
2^{jd} \epsilon\sqrt{\log1/\epsilon} \bigr)
\end{eqnarray*}
but
\[
P \bigl( |\widehat{\beta}_{j,\xi} -\beta_{j,\xi} |\|\psi_{j,\xi}\|
_\infty> \kappa2^{jd} \epsilon\log1/\epsilon \bigr)
\leq \mathrm{e}^{-\fracc{(\kappa2^{jd} \epsilon\sqrt{\log1/\epsilon
})^2}{2(\epsilon2^{j(d-1)/2} \|\psi_{j,\xi}\|_\infty)^2} }
\leq\epsilon^{\kappa^2/2\tau_\infty^2}.
\]
So
\begin{eqnarray*}&&\sum_{j\leq J}\bE \biggl(
\sup_{\xi\in\chi_j}  \bigl( |\beta_{j,\xi}| \| \psi_{j,\xi}\|
_\infty
\Charac{|\widehat{\beta}_{j,\xi} -\beta_{j,\xi} |\|\psi_{j,\xi}\|
_\infty> \kappa2^{jd} \epsilon\sqrt{\log1/\epsilon} }
 \bigr) \biggr)\\
 && \quad \leq C \| f \|_{B^s_{\pi,r} }J \epsilon^{\kappa
^2/2\tau_\infty^2}
\leq C \| f \|_{B^s_{\pi,r} }  \bigl(\epsilon\sqrt{\log1/\epsilon
} \bigr)^{\frace{s-(d+1)/\pi}{s+d-(d+1)/\pi}}
\end{eqnarray*}
if $\kappa^2/2\geq\tau_\infty^2.$
This proves that this term will be of the right order.
Concerning the second term of the second inequality,
\begin{eqnarray*}
&&\sup_{\xi\in\chi_j}  \bigl(|\beta_{j,\xi}| \| \psi_{j,\xi}\|
_\infty
\Charac{| \beta_{j,\xi} |\|\psi_{j,\xi}\|_\infty< 2\kappa2^{jd}
\epsilon\sqrt{\log1/\epsilon} } \bigr)
\\
&& \quad \leq 2\kappa2^{jd} \epsilon\sqrt{\log1/\epsilon} \wedge
C \| f \|_{B^s_{\pi,r} }
2^{-j(s-(d+1)/\pi)},
\end{eqnarray*}
let us again take
\[
2^{j_1} \sim  \bigl(\epsilon\sqrt{\log1/\epsilon} \bigr)^{-\fracc
1{s+d-(d+1)\pi}}
\]
\begin{eqnarray*}&&\sum_{j\leq J} \sup_{\xi\in\chi_j}  \bigl(|\beta
_{j,\xi}| \| \psi_{j,\xi}\|_\infty
\Charac{| \beta_{j,\xi} |\|\psi_{j,\xi}\|_\infty< 2\kappa2^{jd}
\epsilon\log1/\epsilon}  \bigr)
\\ && \quad  \leq
C \| f \|_{B^s_{\pi,r} }  \biggl( \epsilon\log1/\epsilon \sum
_{j\leq j_1} 2^{jd} + \sum_{j_1< j \leq J}
2^{-j(s-(d+1)/\pi)}  \biggr)\\
 && \quad \lesssim
 C \| f \|_{B^s_{\pi,r} }
 \bigl(\epsilon\sqrt{\log
1/\epsilon} \bigr)^{\frace{s-(d+1)/\pi}{s+d-(d+1)/\pi}}.
\end{eqnarray*}
This ends the proof of Theorem \ref{mMinfty}.

\subsection{\texorpdfstring{Proof of  Theorem \protect\ref{mM}}{Proof of  Theorem 1}}\label{sec7.2}

As in the previous proof, we begin with the decomposition,
\[
\|\hat{f_\epsilon}-f\|^p_p \leq 2^{p-1}  \bigl( \|\hat{f_\epsilon}-
A_J(f)\|^p_p +\|A_J(f) -f\|^p_p  \bigr).
\]

Using the fact that, for $ \pi\geq p$, $ B^s_{\pi,r} \subset
B^s_{p,\infty},$ and for
$ \pi\leq p ,  B^s_{\pi,r} \subset B^{s-(d+1)(1/\pi-
1/p)}_{p,\infty}$, we obtain
\[
\|A_J(f) -f\|_p^p \leq C \|f\|^p_{B^s_{\pi,r}} 2^{-Jsp} \qquad \mbox{if }
 \pi\geq p
\]
and
\[
\|A_J(f) -f\|_p^p \leq C \|f\|^p_{B^s_{\pi,r}} 2^{-J(s-(d+1)(1/\pi
-1/p)) p}\qquad \mbox{if }  \pi\leq p.
\]
We have
$2^{-Jsp}
\le(\epsilon\sqrt{\log1/\epsilon})^{\fracc
{sp}{d-\fraca12}} \le(\epsilon\sqrt{\log1/\epsilon})^{\fracc
{sp}{s+d-\fraca12}} $.
Obviously, this term has the right rate for $\pi\geq p$.
For $\pi<p$,
\begin{eqnarray*}
2^{-J(s-(d+1)(1/\pi -1/p)) }&\le&(\epsilon\log1/\epsilon)^{\fracc
{(s-(d+1)(1/\pi -1/p))}{d-\fraca12}}\\
&\le&(\epsilon\log1/\epsilon
)^{\fracc
{(s-(d+1)(1/\pi -1/p))}{s+d-(d+1)/\pi}},
\end{eqnarray*}
thanks to $s\ge(d+1)/\pi-1/2$. This gives the right rate for $dp>d+1.$
For $dp\le d+1$, we have (again as $s\ge(d+1)/\pi-1/2$ ),
$s-(d+1)(\frac1\pi-\frac1p)\ge d-\frac12$, so
$2^{-J(s-(d+1)(1/\pi -1/p)) }\le(\epsilon\log1/\epsilon)^{\fracc
{(s-(d+1)(1/\pi -1/p))}{d-\fraca12}}\le(\epsilon\log1/\epsilon
)^{\fracc
{s}{s+d-\fraca12}}$. Finally, this proves that the bias term above
always has the right rate.

Let us now investigate the stochastic term
\[
\bE\|\hat{f}- A_J(f)\|^p_p \leq C J^{p-1}\sum_{j<J} \bE\biggl\| \sum
_{\xi\in\chi_j}  \bigl( \widehat{\beta}_{j,\xi}
\Charac{|\widehat{\beta}_{j,\xi} | \geq\kappa2^{j\nu} \epsilon\sqrt
{\log1/\epsilon}}
-\beta_{j,\xi}  \bigr)\psi_{j,\xi} \biggr\|^p_p.
\]
But
\begin{eqnarray*}&&\biggl\| \sum_{\xi\in\chi_j}  \bigl ( \widehat{\beta}_{j,\xi}
\Charac{|\widehat{\beta}_{j,\xi} | \geq\kappa2^{j\nu} \epsilon\sqrt
{\log1/\epsilon}}
-\beta_{j,\xi}  \bigr)\psi_{j,\xi} \biggr\|^p_p
\\
&& \quad \leq C 2^{p-1}  \biggl( \biggl\| \sum_{\xi\in\chi_j}  \bigl( (\widehat{\beta
}_{j,\xi} -\beta_{j,\xi} )\psi_{j,\xi}
\Charac{|\widehat{\beta}_{j,\xi} |\geq\kappa2^{j\nu} \epsilon\sqrt
{\log1/\epsilon} }
 \bigr) \biggr\|^p_p\\
&&\hphantom{ C 2^{p-1}  \biggl(} \qquad {}  + \biggl\| \sum_{\xi\in\chi_j} \bigl ( \beta_{j,\xi}
\psi_{j,\xi}
\Charac{|\widehat{\beta}_{j,\xi} | < \kappa2^{j\nu} \epsilon\sqrt
{\log1/\epsilon}}
 \bigr) \biggr\|^p_p  \biggr)
\\ && \quad \leq C
 \biggl( \sum_{\xi\in\chi_j} |\widehat{\beta}_{j,\xi} -\beta_{j,\xi
} |^p \|\psi_{j,\xi}\|^p_p
\Charac{|\widehat{\beta}_{j,\xi}| \geq\kappa2^{j\nu} \epsilon
\sqrt{\log1/\epsilon} } \\
&&\hphantom{C
 \biggl(} \qquad {}
+ \sum_{\xi\in\chi_j} |\beta_{j,\xi}|^p \| \psi_{j,\xi}\|^p_p
\Charac{|\widehat{\beta}_{j,\xi} < \kappa2^{j\nu} \epsilon\sqrt
{\log1/\epsilon} }
 \biggr).
\end{eqnarray*}
In turn,
\begin{eqnarray*}
&&|\widehat{\beta}_{j,\xi} -\beta_{j,\xi} |^p
 \Charac{|\widehat{\beta}_{j,\xi} | \geq\kappa2^{j\nu} \epsilon
\sqrt{\log1/\epsilon}}
\\
&& \quad
= |\widehat{\beta}_{j,\xi} -\beta_{j,\xi} |^p
\Charac{|\widehat{\beta}_{j,\xi} | \geq\kappa2^{j\nu} \epsilon\sqrt
{\log1/\epsilon}}
 \bigl( \Charac{|\beta_{j,\xi} |\geq\fracd{\kappa}2 2^{j\nu}
\epsilon\sqrt{\log1/\epsilon} }
+ \Charac{\beta_{j,\xi} | < \fracd{\kappa}2 2^{j\nu} \epsilon
\sqrt{\log1/\epsilon} }  \bigr)
\\
&& \quad
\leq
|\widehat{\beta}_{j,\xi} -\beta_{j,\xi} |^p
\Charac{|\widehat{\beta}_{j,\xi}-\beta_{j,\xi}| \geq\fracd{\kappa}2
2^{j\nu} \epsilon\sqrt{\log1/\epsilon} }
+
|\widehat{\beta}_{j,\xi} -\beta_{j,\xi} |^p
\Charac{| \beta_{j,\xi} | >\fracd{\kappa}2 2^{jd} \epsilon\sqrt
{\log1/\epsilon}}
\end{eqnarray*}
and
\begin{eqnarray*}
|\beta_{j,\xi}|^p
\Charac{|\widehat{\beta}_{j,\xi} < \kappa2^{jd} \epsilon\sqrt{\log
1/\epsilon} }
&=&
|\beta_{j,\xi}|^p
\Charac{|\widehat{\beta}_{j,\xi} < \kappa2^{jd} \epsilon\sqrt{\log
1/\epsilon} } \\
&&{}  \times
 \bigl( \Charac{ | \beta_{j,\xi} | \geq2 \kappa2^{jd} \epsilon
\sqrt{\log1/\epsilon} }
+\Charac{ | \beta_{j,\xi} | < \kappa2^{jd} \epsilon\sqrt{\log
1/\epsilon} }
 \bigr) \\&\leq&
|\beta_{j,\xi}|^p
\Charac{|\widehat{\beta}_{j,\xi} -\beta_{j,\xi} | > \kappa2^{jd}
\epsilon\sqrt{\log1/\epsilon} }
+ |\beta_{j,\xi}|^p
\Charac{| \beta_{j,\xi} | < \kappa2^{jd} \epsilon\sqrt{\log
1/\epsilon} }.
\end{eqnarray*}

We now have the following bound
using direct computation:
\[
\bE \bigl(|\widehat{\beta}_{j,\xi} -\beta_{j,\xi} |^p
\Charac{|\widehat{\beta}_{j,\xi}-\beta_{j,\xi} | \geq\fracd{\kappa
}2 2^{j\nu} \epsilon\sqrt{\log1/\epsilon} }  \bigr)
\leq C 2^{jp\nu} \epsilon^p  \bigl(\kappa \sqrt{\log1/\epsilon
} \bigr)^{p-1} \epsilon^{\kappa^2/2}.
\]
Hence,
\begin{eqnarray*}
&&J^{p-1} \sum_{j\leq J}\sum_{\xi\in\chi_j} \bE \bigl( |\widehat{\beta
}_{j,\xi} -\beta_{j,\xi} |^p  \|\psi_{j,\xi}\|_p^p
\Charac{|\widehat{\beta}_{j,\xi}-\beta_{j,\xi}| \geq\fracd{\kappa}2
2^{j\nu} \epsilon\sqrt{\log1/\epsilon} }  \bigr)
\\
&& \quad \leq C J^{p-1} \sum_{j\leq J} \epsilon^p  \bigl(\kappa \sqrt{\log
1/\epsilon} \bigr)^{p-1} \epsilon^{\kappa^2/2}
\sum_{\xi\in\chi_j} \|\psi_{j,\xi}\|_p^p\\
&& \quad \leq C \sqrt{\log
1/\epsilon}^{p-1} \epsilon^p  \bigl(\kappa \sqrt{\log1/\epsilon
} \bigr)^{p-1} \epsilon^{\kappa^2/2} \sum_{j\leq J}
2^{jp\nu} 2^{j (dp/2 + (p/2-2)_+)} \\
&& \quad \leq C\epsilon^p
\end{eqnarray*}
if $\kappa\ge\sqrt{2p}$ is large enough (we have used \eqref{norm-est}).

Using the forthcoming inequality \eqref{FOND2},
\begin{eqnarray*}&& \sum_{\xi\in\chi_j} \bE \bigl( |\widehat{\beta
}_{j,\xi} -\beta_{j,\xi} |^p   \|\psi_{j,\xi}\|_p^p
\Charac{| \beta_{j,\xi} | > \fracd{\kappa}2 2^{j\nu} \epsilon
\sqrt{\log1/\epsilon} }
 \bigr)\\
&& \quad \le C \epsilon^p \sum_{\xi\in\chi_j} 2^{j\nu p} \|\psi_{j,\xi}\|_p^p
\Charac{| \beta_{j,\xi} | > \fracd{\kappa}2 2^{j\nu} \epsilon
\sqrt{\log1/\epsilon} }
\\
&& \quad \leq C \epsilon^p  \biggl( \frac{\kappa}2 \epsilon\sqrt{\log
1/\epsilon}  \biggr)^{-q}.
\end{eqnarray*}

Hence,
\begin{eqnarray*} &&J^{p-1}\sum_{j\leq J } \sum_{\xi\in\chi_j} \bE
\bigl( |\widehat{\beta}_{j,\xi} -\beta_{j,\xi} |^p  \|\psi_{j,\xi}\|_p^p
\Charac{| \beta_{j,\xi} |> \fracd{\kappa}2 2^{j\nu} \epsilon\sqrt
{\log1/\epsilon} }
 \bigr)
\\ && \quad \leq C J^p \sqrt{\log1/\epsilon}^{-p}  \biggl( \frac{\kappa}2
 \biggr)^{-q}  \bigl( \epsilon\sqrt{\log1/\epsilon}  \bigr)^{p-q}
\end{eqnarray*}
and, as $J\leq\frac{1}{d-1/2}\log1/\epsilon+1$,
\[ \leq C  \bigl( \epsilon\sqrt{\log1/\epsilon}  \bigr)^{p-q}
(\log1/\epsilon)^{p/2}.
\]
This term is thus of the right rate. Let us now turn to
\begin{eqnarray*} &&\sum_{\xi\in\chi_j} \bE \bigl( |\beta_{j,\xi}|^p
\|  \psi_{j,\xi}\|_p^p
\Charac{|\widehat{\beta}_{j,\xi} -\beta_{j,\xi} |> \kappa2^{j\nu}
\epsilon\sqrt{\log1/\epsilon} }
 \bigr)\\ && \quad = \sum_{\xi\in\chi_j} |\beta_{j,\xi}|^p \| \psi
_{j,\xi}\|_p^p P \bigl( |\widehat{\beta}_{j,\xi} -\beta_{j,\xi} |>
\kappa2^{j\nu} \epsilon\sqrt{\log1/\epsilon} \bigr).
\end{eqnarray*}
As the standard deviation of $\widehat{\beta}_{j,\xi}
-\beta_{j,\xi}$ is smaller than
$\eps2^{j(d-1)/2}$,
\[
P \bigl( |\widehat{\beta}_{j,\xi} -\beta_{j,\xi} |> \kappa2^{j\nu}
\epsilon\sqrt{\log1/\epsilon} \bigr)
\leq\epsilon^{\kappa^2/2}.
\]
So
\begin{eqnarray*}&& J^{p-1}\sum_{j\leq J}\bE \biggl(\sum_{\xi\in\chi_j}
  |\beta_{j,\xi}|^p \| \psi_{j,\xi}\|_p^p
\Charac{|\widehat{\beta}_{j,\xi} -\beta_{j,\xi}| > \kappa2^{j\nu}
\epsilon\sqrt{\log1/\epsilon} }
 \biggr)\\
&& \quad \leq C \| f \|_{B^s_{\pi,r} }^p J^{p-1} \epsilon^{\kappa^2/2}
\leq C  \bigl(\epsilon\sqrt{\log1/\epsilon} \bigr)^p
\end{eqnarray*}
if $\kappa^2$ is large enough
(where we have used that $B^s_{\pi,r} \subset B^{s'}_{p,r} \subset
\bL_p$ with $s'=s-(d+1)(1/\pi-1/p)$).
Hence this term also is of the right order.

Let us turn now to the last one, using (\ref{FOND2})
\[
\sum_{\xi\in\chi_j} |\beta_{j,\xi}|^p \| \psi_{j,\xi}\|_p^p
\Charac{| \beta_{j,\xi} | < 2\kappa2^{j\nu} \epsilon\sqrt{\log
1/\epsilon} }
\leq  \bigl( 2\kappa\epsilon\sqrt{\log1/\epsilon}  \bigr)^{p-q}.
\]
Hence,
\begin{eqnarray*}&& J^{p-1}\sum_{j\leq J} \sup_{\xi\in\chi_j}  \bigl(
|\beta_{j,\xi}|  \|^p \psi_{j,\xi}\|_p^p
\Charac{| \beta_{j,\xi} |< 2\kappa2^{j\nu} \epsilon\sqrt{\log
1/\epsilon} }  \bigr)
\\
&& \quad
\leq C \| f \|_{B^s_{\pi,r} } J^p  \bigl( \epsilon\sqrt{\log
1/\epsilon}  \bigr)^{p-q}
\leq C \| f \|_{B^s_{\pi,r} } \sqrt{\log1/\epsilon}^p  \bigl(
\epsilon\sqrt{\log1/\epsilon}  \bigr)^{p-q}.
\end{eqnarray*}
This proves that all the terms have the proper rate.
It remains now to state and prove the following lemma.
\begin{lemma}\label{FOND}
Let
$ A= \{ (s,\pi),   s > (d+1)(\frac1\pi- \frac1p) \cap(s>0)
 \},
$
and $f \in B^s_{\pi,r},   1\leq\pi\leq\infty,  1\leq p
<\infty.$
If
$ \sum_{\xi\in\chi_j}( |\beta_{j,\xi}| \| \psi_{j,\xi}\|_\pi
)^\pi=
\rho_j^\pi2^{-js \pi}$ with $\rho\in l_r(\bN)$,
then, with $\nu= \frac{d-1}2$,
\[
\sum_{\xi\in\chi_j} \biggl(\frac{ |\beta_{j,\xi}|}{2^{j\nu}}
 \biggr)^q
 (2^{j\nu}\| \psi_{j,\xi}\|_p )^p \leq C \rho_j^{q},
\]
where $q < p$ is as follows:
\begin{longlist}[(3)]
\item[(1)]
$ p-q = \frac{ sp}{s+d-1/2}  (q=\frac{(d-1/2)p}{s+d-1/2} ) $ in
the following domain I:
\[
 \bigl\{ (s,\pi),  \bigl( s(1/p-1/4) \geq (d-1/2) (1/\pi-1/p) \bigr) \cap
A \bigr\}.
\]
Moreover, we have the following slight modification at the frontier:
the domain becomes
\[
 \bigl\{ (s,\pi),  \bigl( s(1/p-1/4) =(d-1/2)(1/\pi-1/p)\bigr) \cap A
\bigr\}
\]
and the inequality
\[
\sum_{\xi\in\chi_j} \biggl(\frac{ |\beta_{j,\xi}|}{2^{j\nu}}
 \biggr)^q  (2^{j\nu}\| \psi_{j,\xi}\|_p )^p\leq
C \rho_j^q j^{1-q/\pi}.
\]

\item[(2)]
$ p-q = \frac{(s-2(1/\pi-1/p))p}{s +d -2/\pi}$  $(q=\frac
{dp+2}{s+d-2/\pi} ) $ in the following domain II:
\[
 \bigl\{ (s,\pi)  \bigl(s > dp(1/\pi-1/p)\bigr) \cap\bigl(s(1/p-1/4) <
(d-1/2)(1/\pi-1/p)\bigr) \cap A \bigr\}.
\]
\item[(3)]
$ p-q = \frac{(s- (d+1)(1/\pi-1/p))p}{s +d -(d+1)/\pi} $  $(q=\frac
{dp-(d+1)}{s+d-(d+1)/\pi} ) $ in the following domain III:
\[
\biggl\{ (s,\pi),  \biggl( dp\biggl(\frac1\pi-\frac1p\biggr) \geq s \biggr)\cap A ,
\mbox{ for }  \frac1p < \frac d{d+1}\biggr\}.
\]
\end{longlist}
\end{lemma}

This lemma is to be used essentially through the following
corollary.\vadjust{\goodbreak}
\begin{corollary}
Respectively, in the domains I, II, III, we have, for $q$ described in
the lemma and $f \in B^s_{\pi,r}$,
%
\begin{equation}\label{FOND1}
\sum_{\xi\in\chi_j} \Charac{\fraca{ |\beta_{j,\xi}| }{ 2^{j\nu}
} \geq\lambda}( 2^{j\nu}\|\psi_{j,\xi}\|_p)^p
\leq C \rho_j^q \lambda^{-q},
\end{equation}
%
\begin{equation}\label{FOND2}
\sum_{ \xi\in\chi_j } \Charac{ \fraca{ |\beta_{j, \xi} | }{
2^{j\nu} } \leq2^{j\nu} \lambda}
|\beta_{j,\xi}|^p \| \psi_{j,\xi}\|_p^p \leq C \rho_j^q \lambda^{p-q}
\end{equation}
with an obvious modification for
\[
\bigl \{ (s,\pi), \bigl ( s(1/p-1/4) =(d-1/2)(1/\pi-1/p)\bigr) \cap A
\bigr\}.
\]
\end{corollary}

\begin{pf}
Let us recall that on a measure space $(X,\mu)$ we have, if $h \in\bL
_q(\mu)$ then
$\mu(|h| \geq\lambda) \leq\frac{\|h\|_q^q}{\lambda^q}$
and, as $q<p,$
\begin{eqnarray*}
\int_{|h| \leq\lambda} |h|^p \,\mathrm{d}\mu&\leq&\int (|h| \wedge\lambda)^p
\,\mathrm{d}\mu=\int_{0}^\lambda
px^{p-1}\mu( |h| \geq x)\,\mathrm{d}x\\
 &\leq&\int_{0}^\lambda
px^{p-1} \frac{\|h\|_q^q}{x^q}\,\mathrm{d}x = \frac{p \|h\|_q^q}{p-q}
\lambda^{p-q}.
\end{eqnarray*}
For the corollary, we take $X=\chi_j$, $\mu(\xi) =(2^{j\nu}\|\psi
_{j,\xi}\|_p)^p
$
and $h(\xi) = \frac{ |\beta_{j,\xi}| }{ 2^{j\nu}} .$
\end{pf}

\begin{pf*}{Proof of Lemma \ref{FOND}}
Let us fix $q$ (chosen later) and investigate separately the two cases
$ q\geq\pi$ and $q<\pi.$

For $q\geq\pi$, we have, using \eqref{norm-est},
\begin{eqnarray*}
I_j(f, q,p)&=&\sum_{\xi\in\chi_j}  \biggl| \frac{\beta_{j,\xi
}}{2^{j\nu}} \biggr|^q \| 2^{j\nu} \psi_{j,\xi}\|_p^p
\sim
2^{j \nu(p-q) } \sum_{\xi\in\chi_j } | \beta_{j,\xi}|^q
\biggl(\frac{2^{jd}}{W_j(\xi) } \biggr)^{p/2-1}\\
&\leq&
2^{j \nu(p-q) }
 \biggl( \sum_{\xi\in\chi_j } \biggl( | \beta_{j,\xi}|^q
\biggl(\frac{2^{jd}}{W_j(\xi)} \biggr)^{p/2-1} \biggr)^{\pi/q}
\biggr)^{q/\pi}\\
&=&
2^{j \nu(p-q) }
 \biggl(
\sum_{\xi\in\chi_j }| \beta_{j,\xi}|^\pi  \biggl(\frac
{2^{jd}}{W_j(\xi) } \biggr)^{(p/2-1)\pi/q }  \biggr)^{q/\pi}\\
&=&
2^{j \nu(p-q) }
 \biggl( \sum_{\xi\in\chi_j }| \beta_{j,\xi}|^\pi  \biggl(\frac
{2^{jd}}{W_j(\xi) } \biggr)^{\pi/2-1}
 \biggl(\frac{2^{jd}}{W_j(\xi) } \biggr)^{ \fracd\pi q(p/2-1) -(\pi
/2-1) }  \biggr)^{q/\pi}\\
&\leq&
2^{j \nu(p-q) }
 \biggl( \sum_{\xi\in\chi_j }| \beta_{j,\xi}|^\pi  \biggl(\frac
{2^{jd}}{W_j(\xi) } \biggr)^{\pi/2-1}
 \biggr)^{q/\pi} 2^{j(d+1) (\fracb{p-q}2+ q(\fraca1\pi- \fraca1q)}
\end{eqnarray*}
if we choose $q$ such that
$(sq + d(p-q) + (d+1)q(\frac1\pi- \frac1q) ) =0$ gives $q= \frac
{pd-(d+1)}{s+d-(d+1)\pi} $. Hence,
$p-q = \frac{s-(d+1)(1/\pi- 1/p)}{s+d-(d+1)\pi}p $; $ q-\pi= -
\pi\frac{s-pd(1/\pi- 1/p)}{s+d-(d+1)\pi}.$
Thus, in domain~III,
\[
 \biggl\{ \biggl(\frac1p < \frac d{d+1}\biggr) \cap \bigl(s-(d+1)(1/\pi- 1/p) >0\bigr) \cap
\bigl(s-pd(1/\pi- 1/p) \leq0\bigr)  \biggr\},
\]
we have $0<q<p$, $\pi\leq q$,  $\sum_{\xi\in\chi_j}  |
\frac{\beta_{j,\xi}}{2^{j\nu}} |^q \| 2^{j\nu} \psi_{j,\xi
}\|_p^p \leq\rho_j^q$.

For $q< \pi$, we have, using \eqref{norm-est},
\[
I_j(f, q,p)=\sum_{\xi\in\chi_j}  \biggl| \frac{\beta_{j,\xi
}}{2^{j\nu}} \biggr|^q \| 2^{j\nu} \psi_{j,\xi}\|_p^p
\sim
2^{j \nu(p-q) } \sum_{\xi\in\chi_j } | \beta_{j,\xi}|^q
 \biggl(\frac{2^{jd}}{W_j(\xi) } \biggr)^{p/2-1},
\]
\begin{eqnarray*} &&2^{j \nu(p-q) }  \sum_{\xi\in\chi_j } | \beta
_{j,\xi}|^q  \biggl(\frac{2^{jd}}{W_j(\xi) } \biggr)^{(\pi/2-1)q/\pi}
 \biggl(\frac{2^{jd}}{W_j(\xi) } \biggr)^{ (p/2-1)-(\pi/2-1)q/\pi}\\
&& \quad \leq
2^{j \nu(p-q) }
\biggl ( \sum_{\xi\in\chi_j }| \beta_{j,\xi}|^\pi  \biggl(\frac
{2^{jd}}{W_j(\xi) } \biggr)^{\pi/2-1}
 \biggr)^{q/\pi}\\
 && \qquad {}\times  \biggl( \sum_{\xi\in\chi_j } \biggl (\frac
{2^{jd}}{W_j(\xi) } \biggr)^{ \fracc\pi{\pi-q} ((p/2-1)-(\pi
/2-1)q/\pi)} \biggr)^{1-q/\pi}\\
&& \quad  \sim 2^{j \nu(p-q) }
 \biggl( \sum_{\xi\in\chi_j } | \beta_{j,\xi} |^\pi \| \psi
_{j,\xi} \|_\pi^\pi
 \biggr)^{q/\pi}  \biggl( \sum_{\xi\in\chi_j } \biggl(\frac
{2^{jd}}{W_j(\xi) }\biggr)^{ \fracd{\pi(p-q)}{(2(\pi-q))}-1}
\biggr)^{1-q/\pi}
\\
&& \quad  \sim 2^{j \nu(p-q) }
 \biggl( \sum_{\xi\in\chi_j } | \beta_{j,\xi} |^\pi \| \psi
_{j,\xi} \|_\pi^\pi
 \biggr)^{q/\pi}  \biggl( \sum_{\xi\in\chi_j } \| \psi_{j,\xi} \|
_{\fraca{\pi(p-q)}{(\pi-q)}}^{\fraca{\pi(p-q)}{(\pi-q)}}
 \biggr)^{1-q/\pi}.
\end{eqnarray*}

Now let us investigate separately the cases $
\frac{\pi(p-q)}{(\pi-q)}$ smaller, greater or equal to $4$.

\textit{Case}
$ \frac{\pi(p-q)}{(\pi-q)} <4$.
Using \eqref{norm-est}--\eqref{BI}, we have
\[
I_j(f, q,p)\leq C 2^{j \nu(p-q) } \rho_j^q 2^{-jsq}
2^{jd(p-q)/2} \leq C \rho_j^q.
\]
If we define $q$ such that $ -sq + (p-q)(d-1/2)=0$, that is, $q=\frac
{p(d-1/2)}{s+d-1/2}$, then
$p-q= \frac{sp}{s+d-1/2} >0$.
So
$ \pi-q =\pi\frac{s -(d-1/2)p(1/\pi-1/p) }{s+d-1/2} >0
\Leftrightarrow \frac sp >(d-1/2)(1/\pi-1/p) $
and
$ \frac{\pi(p-q)}{(\pi-q)} =\frac{sp}{s -(d-1/2)p(1/\pi-1/p)} <4
\Leftrightarrow s(1/p-1/4) > (d-1/2)(1/\pi-1/p) $.
Hence we only need to impose $s(1/p-1/4) > (d-1/2)(1/\pi-1/p)$ and
domain I is given by
\[
 \bigl\{ \bigl(s-(d+1)(1/\pi- 1/p) >0\bigr) \cap(s>0) \bigr\} \cap
\{ s(1/p-1/4) > (d-1/2)(1/\pi-1/p)  \}
\]
on which
$ I_j(f, q,p)\leq C \rho_j^{\fracc{p(d-1/2)}{s+d-1/2} }$.

\textit{Case}
$ \frac{\pi(p-q)}{(\pi-q)} >4$.
Using \eqref{norm-est}--\eqref{BI}, we have
\[
I_j(f, q,p)\leq C 2^{j \nu(p-q) } \rho_j^q 2^{-jsq} 2^{jd(p-q)/2}
2^{j(\fraca{(p-q)}{2} -2\fracb{\pi-q}\pi)}.
\]
If we put
$(p-q)d -sq -2\frac{\pi-q}\pi=0 \Leftrightarrow q=\frac
{pd-2}{s+d-2/\pi}$,
we have
$p-q =\frac{s-2(1/\pi- 1/p)}{s+d-2/\pi}p >0 \Leftrightarrow
s-2(1/\pi- 1/p) >0$ and
$\pi-q =\frac{s-dp(1/\pi- 1/p)}{s+d-2/\pi}\pi >0 \Leftrightarrow
s-dp(1/\pi- 1/p) >0$.

Moreover,
$ \frac{\pi(p-q)}{(\pi-q)} = \frac{s-2(1/\pi- 1/p)}{s-dp(1/\pi
-1/p)}p >4 \Leftrightarrow s(1/p-1/4) < (d-1/2)(1/\pi-1/p) $.
Hence, on the domain
\begin{eqnarray*}
 &&\bigl\{ (0<s)\cap \bigl( s> (d+1)(1/\pi-1/p)\bigr)\\
 &&\hphantom{\bigl\{} \cap\bigl(s > dp(1/\pi-1/p)\bigr)
\cap\bigl(s(1/p-1/4) < (d-1/2)(1/\pi-1/p)\bigr) \bigr\},
\end{eqnarray*}
we have
$I_j(f,q,q) \leq C \rho_j^{\frace{pd-2}{s+d-2/\pi}}$.

\textit{Case}
$ \frac{\pi(p-q)}{(\pi-q)} =4$.
Using \eqref{norm-est}--\eqref{BI}, we have
\[
I_j(f, q,p)\leq C 2^{j \nu(p-q) } \rho_j^q 2^{-jsq} j^{1-q/\pi} 2^{jd(p-q)/2}
\leq C \rho_j^q j^{1-q/\pi}
\]
if
$(p-q)(d-1/2) -sq =0 \Leftrightarrow q=p \frac{d-1/2}{s+d-1/2} $.
This is realized
either if $p=4=\pi$ and for $s>0$
or if
$p\neq4, \pi\neq4$ and $ 0< q= \pi\frac{p-4}{\pi-4} = p \frac
{d-1/2}{s+d-1/2}.$
Moreover,
$q<\pi$ and $q<p$ $\Leftrightarrow 4 <p<\pi;  \mbox{ or }
4< \pi<p$
and
$ \frac{\pi(p-q)}{(\pi-q)} =4 \Leftrightarrow s(1/p-1/4)
=(d-1/2)(1/\pi-1/p) \}$.
Hence,
on the domain
\[
 \bigl\{ (s,\pi),  \bigl( s(1/p-1/4) =(d-1/2)(1/\pi-1/p)\bigr) \cap(s>0)
\cap
\bigl(s> (d+1)(1/\pi- 1/p)\bigr)  \bigr\},
\]
we have
$ I_j(f, q,p)\leq C \rho_j^{\frace{p(d-1/2)}{s+d-1/2}} j^{\fracc{s}{s+d-1/2}}$.
\end{pf*}

\section{Proof of the lower bounds}\label{sec8}
In this section, we prove the lower bounds:
for $ 0<s<\infty,  1\leq\pi\leq\infty,  0<r\leq\infty,
 0<M<\infty$, denoting by
$ B^{s}_{\pi,r}(M)$ the ball of radius $M$ of the space $B^{s}_{\pi
,r}$ and, by $\cEstim$, the set of all estimators, we consider
\begin{eqnarray*}
\omega_p(s,\pi,r, M,\epsilon) &=&\inf_{f^\star\in\cEstim} \sup
_{f \in B^{s}_{\pi,r}(M)} \bE\|f^\star-f\|_p^p,
\\
\omega_\infty(s,\pi,r, M,\epsilon) &=&\inf_{f^\star\in\cEstim}
\sup_{f \in B^{s}_{\pi,r}(M)} \bE\|f^\star-f\|_\infty.
\end{eqnarray*}
The main tool will be the classical lemma introduced by Fano in 1952
\cite{fano}.
We will use the version of Fano's lemma introduced in \cite{birge}.
For details on general lower bound results, see also \cite
{tsybakov08introdnonparestim}.
Let us recall that $K(P,Q)$ denotes the Kullback information
``distance''
between $P$ and $Q$.

\begin{lemma}
Let $\cA$ be a sigma algebra on the space $\Omega, $ and
$A_i \in\cA$,   $i\in\{ 0,1,\ldots,m \}$, such that $\forall i\neq
j $,   $A_i\cap A_j =\emptyset,$
$P_i,   i\in\{0,1,\ldots,m\}$ are $m+1$ probability measures on
$(\Omega,\cA).$ Define
\[
p := \sup_{i=0}^m P_{i}( A^c_i)
\quad\mbox{and}\quad\kappa
:= \inf_{j \in\{0,1,\ldots ,m\}}\frac1m \sum_{i\neq j} K(P_i,P_j),
\]
then
%
\begin{equation}
p \geq\frac12 \wedge \bigl (C \sqrt m \exp(- \kappa) \bigr)
\qquad\mbox{with }  C= \exp \biggl(-\frac3e \biggr).
\end{equation}
\end{lemma}

This inequality will be used in the following way:
Let $H_\epsilon$ be the Hilbert space of measurable functions on $Z=
\bS^{d-1}\times[-1,1]$ with the
scalar product
\[
\langle\phi,\psi\rangle_\epsilon= \epsilon^2 \int_{\bS^{d-1}}
\int_{-1}^1 \phi(\theta, s) \psi(\theta,s)\,\mathrm{d}\sigma(\theta)
\frac{\mathrm{d}s}{(1-s^2)^{(d-1)/2}}.
\]

It is well known
that there
exists a (unique) probability measure on $(\Omega,\cA) \dvt Q_f $ the
density of which, with respect to $P,$ is
\[
\frac{\mathrm{d}Q_f}{\mathrm{d}P} = \exp\biggl(W^\epsilon(f) - \frac12\| f \|^2_{H_\epsilon
}\biggr) .
\]
Let us now choose $f_0, f_1,\ldots,f_m$ in $B^{s}_{\pi,r}(M)$ such that
$ i\neq j \Longrightarrow\| f_i- f_j\|_p \geq\delta$ and denote $P_i
=Q_{R(\fraca{f_i}{\epsilon^2})}$. Let $f^\star$ be an arbitrary
estimator of $f$.
Obviously, the sets $A_i = (\| f^\star-f\|_p < \frac\delta2)$ are
disjoint sets
and we have, for $i\neq j$,
\[
K(P_i, P_j)
= \frac1{2\epsilon^2} \int_{Z} |R(f_i-f_j) |^2\,\mathrm{d}\mu.
\]
Now
\begin{eqnarray*}
\omega_p(s,\pi,q, M,\epsilon) &\geq&
\inf_{f^\star\in\cEstim} \sup_{f_i, i=0,1,\ldots,m} \bE\|f^\star-f_i\|
_p^p\\
&\geq& \biggl(\frac\delta2 \biggr)^p
\inf_{f^\star\in\cEstim} \sup_{f_i, i=0,1,\ldots,m} P \biggl(\|f^\star
-f_i\|_p
\geq\frac\delta2 \biggr).
\end{eqnarray*}
Likewise,
\[
\omega_\infty(s,\pi,q, M,\epsilon) \geq \biggl(\frac\delta2 \biggr)
\inf_{f^\star\in\cEstim} \sup_{f_i, i=0,1,\ldots,m} P\biggl (\|f^\star
-f_i\|_\infty\geq\frac\delta2 \biggr).
\]
Using Fano's lemma,
\[
\sup_{f_i, i=0,1,\ldots,m} P \biggl(\|f^\star-f_i\|_p \geq\frac\delta
2 \biggr) \geq
\frac12 \wedge  \bigl(C \sqrt m \exp(- \kappa) \bigr)
\]
with
\[
\kappa= \inf_{j=0,\ldots,M}\frac1m \sum_{i\neq j}\frac1{2\epsilon^2}
\int_{Z} |R(f_i-f_j) |^2\,\mathrm{d}\mu.\vspace*{-3pt}
\]

So if, for a given $\epsilon$, we can find $f_0, f_1, \ldots,f_m$ in
$B^{s,0}_{\pi,r}(M)$ such that
$ i\neq j \Longrightarrow\| f_i- f_j\|_p \geq\delta(\epsilon) $ and
$C \sqrt m \exp(- \kappa) \ge1/2,$ then we have
\[
\mbox{for } p<\infty, \qquad \omega_p(s,\pi,q, M,\epsilon)\geq
\tfrac12 \delta(\epsilon)^p    \quad \mbox{and} \quad\omega_\infty(s,\pi,q,
M,\epsilon)   \geq\tfrac12 \delta(\epsilon).\vspace*{-3pt}
\]

In the sequel, we will choose, as usual, sets of functions containing
either two items (sparse case) or a number of order $ 2^{jd}$ or $
2^{j(d-1)}$ (dense cases). We will consider sets of functions that
are basically linear combinations of needlets at a fixed level $f
=\sum_{\xi\in\chi_j} \beta_{j,\xi} \psi_{j,\xi}$. Because the
needlets have different orders of norms, depending on whether they are
around the north pole or closer to the equator, we will have to
investigate different cases. These differences will precisely yield the
different minimax rates.

\vspace*{-3pt}\subsection{Reverse inequality}\vspace*{-3pt}\label{sec8.1}

Because the needlets are not forming an orthonormal system, we cannot
pretend that inequality \eqref{2} is an equivalence. Since, precisely
in the lower-bound evaluations, we need to bound both sides of the $\bL
_p$ norm for terms of the form $ \sum_{\xi\in A_j} \lambda_\xi\psi
_{j,\xi}$ with $A_j \subset\chi_j$. The following section is devoted
to this problem.\vspace*{-3pt}
\begin{proposition}\label{propcomp}
For $A_j \subset\chi_j$,
\[
\frac1{ C } \biggl(\sum_{\xi'\in A_j} \biggl| \biggl\langle\sum_{\xi\in A_j}
\lambda_\xi\psi_{j,\xi} , \psi_{j,\xi'} \biggr\rangle\biggr|^p \| \psi
_{j,\xi'}\|^p_p \biggr)^{1/p} \leq
\biggl\| \sum_{\xi\in A_j} \lambda_\xi\psi_{j,\xi} \biggr\|_p \leq C
 \biggl(\sum_{\xi\in A_j} | \lambda_\xi|^p \| \psi_{j,\xi} \|
^p_p \biggr)^{1/p}.\vspace*{-3pt}
\]
\end{proposition}

\begin{pf}
Let $f = \sum_{\xi\in A_j} \lambda_\xi\psi_{j,\xi}.$ Clearly, by
(\ref{2}),
\[
\biggl\| \sum_{\xi\in A_j} \lambda_\xi\psi_{j,\xi} \biggr\|_p \leq C
 \biggl(\sum_{\xi\in A_j} | \lambda_\xi|^p \| \psi_{j,\xi} \|
^p_p \biggr)^{1/p},\vspace*{-3pt}
\]
and by (\ref{1}),
\[
 \biggl(\sum_{\xi'\in\chi_j} \biggl| \biggl\langle\sum_{\xi\in A_j} \lambda
_\xi
\psi_{j,\xi} , \psi_{j,\xi'} \biggr\rangle\biggr|^p \|
\psi_{j,\xi'}\|^p_p \biggr)^{1/p} \leq C \biggl\| \sum_{\xi\in A_j}
\lambda_\xi\psi_{j,\xi}\biggr\|_p,\vspace*{-3pt}
\]
so obviously,
\begin{eqnarray*}
&&\frac1{C} \biggl (\sum_{\xi'\in A_j} \biggl| \biggl\langle\sum_{\xi\in A_j}
\lambda_\xi\psi_{j,\xi} , \psi_{j,\xi'} \biggr\rangle\biggr|^p \| \psi
_{j,\xi'}\|^p_p \biggr)^{1/p}\\[-2pt]
&& \quad \leq \biggl\| \sum_{\xi\in A_j} \lambda_\xi\psi_{j,\xi} \biggr\|_p \leq C
\biggl (\sum_{\xi\in A_j} | \lambda_\xi|^p \| \psi_{j,\xi} \|
^p_p \biggr)^{1/p}.
\end{eqnarray*}
\upqed\vspace*{-3pt}
\end{pf}\eject

In the sequel, we will look for subset $A_j$ with equilibrated $\bL_p$
norms, that is, such that there exists
\[
0< D_j, \qquad \mbox{such that }   \forall \xi \in A_j,
\| \psi_{j,\xi} \|_p\sim D_j.
\]
(Here and in the rest of this section, $a_{j,\xi}\sim b_j$ will mean that
there exist two absolute constants $c_1$ and $c_2$ -- which will not
be precised for the sake of simplicity -- such that
$c_1 b_j\le a_{j,\xi}\le c_2b_j,  $ for all considered $\xi$.)
As specified above, $D_j$ may have different forms depending on the
regions. Using \eqref{norm-est}, we have
\[
\| \psi_{j,\xi} \|_p \sim \biggl(\frac{2^{jd}}{2^{-j}
+\sqrt{1-|\xi|^2}} \biggr)^{1/2-1/p}.
\]
For our purpose, let us precise Proposition \ref{prop:CUB} by choosing
the cubature points in the following way:
We choose in the hemisphere $\bS_+^d$ strips $S_k=B(A,(2k+1)\eta
)\setminus B(A,2k\eta)$ with $\eta\sim\frac\pi{2 2^{j+1}}$, $k\in
\{0,\ldots,2^j-1\}$ ($A$ is the north pole).
In each of these strips, we choose a maximal $\eta$-net of points
$\tilde\xi$, whose cardinality is of order $k^{d-1}$. Projecting
these points on the ball, we obtain cubature points $ \xi$ on the ball
with coefficients $\lambdaquadra_{j,\xi}\sim2^{-jd}W_j(\xi)$.
As a consequence, we have in the set $\{x \in\bR^d,  |x|\le\frac
1{\sqrt{2}}\}$, about $ 2^{jd}$ points of cubature for which
\[
D_j\sim \|\psi_{j,\xi} \|_p \sim 2^{jd(1/2-1/p)} .
\]

And in the
corona
$\{(1-2^{-2j}\le|x|\le1\}$, we have about $ 2^{j(d-1)}$ points of
cubature for which
\[
D_j\sim\| \psi_{j,\xi} \|_p \sim 2^{j(d+1)(1/2-1/p)} .
\]

Now let us consider a set $A_j$ of cubature points included in one of
the two sets considered just above (either $\{x \in\bR^d,  |x|\le
\frac1{\sqrt{2}}\}$ or $\{(1-2^{-2j}\le|x|\le1\}$). Consider also
the matrix (parametrized by $A_j$)
\[
\cM(A_j) =  ( \langle\psi_{j,\xi}, \psi_{j,\xi'} \rangle
 )_{\xi, \xi' \in A_j \times A_j}.
\]
We have, for any $\lambda\in l_p(A_j)$, using Proposition \ref{propcomp},
\[
\| \cM(A_j)(\lambda)\|_{l_p(A_j)} \leq C' \| \lambda\|_{l_p(Aj)}.
\]
On the other hand, let us observe that, using \eqref{A},
\[
0< c\leq\| \psi_{j,\xi} \|_2^2 = \langle\psi_{j,\xi},
\psi_{j,\xi} \rangle\leq1.
\]
Thus,
\[
\cM(A_j) = \Diag(\cM(A_j)) + \cM'(A_j) =\Diag(\cM(A_j))\bigl(\mathrm{Id}
+[\Diag(\cM(A_j))]^{-1 } \cM'(A_j) \bigr),
\]
where
$\Diag(\cM(A_j)) $ is the diagonal matrix parametrized by $A_j$
extracted from
$\cM(A_j)$. Clearly, each of the terms of $[\Diag(\cM(A_j)]^{-1 } $
is bounded by $c^{-1}$.

So if $ \| [\Diag(\cM(A_j)]^{-1 } \cM'(A_j) \|_{\cL(l_p(A_j))} \leq
\alpha<1$, we have
\[
\| \cM(A_j)^{-1 } \|_{\cL(l_p(A_j))} \leq c^{-1} \frac
1{1-\alpha}.
\]
Let us prove that we can choose $A_j$ large enough and such that such
an $\alpha$ exists.
Using the Schur lemma (see \cite{grafakos}, Appendix 29),
\[
\| [\Diag(\cM(A_j)]^{-1 } \cM'(A_j) \|_{\cL(l_p(A_j))} \leq
c^{-1}\sup_{\xi\in A_j} \sum_{\xi' \neq\xi, \xi' \in A_j}
| \langle\psi_{j,\xi}, \psi_{j,\xi'} \rangle|.
\]
Now, using \eqref{EqFond2},
\[
| \langle\psi_{j,\xi}, \psi_{j,\xi'} \rangle|\leq C_M^2
\int_{B^d} \frac1{\sqrt{W_j(x)}} \frac1{(1+2^jd(x,\xi))^M}\frac
1{\sqrt{W_j(x)}} \frac1{(1+2^jd(x,\xi'))^M}\,\mathrm{d}x
\]
and thus, by triangular inequality,
\begin{eqnarray*}
| \langle\psi_{j,\xi}, \psi_{j,\xi'} \rangle|&\leq&
\frac{C_M^2}{(1+2^j d(\xi, \xi'))^M} \int_{B^d}\frac1{2^{-j}
+\sqrt{1-|x|^2}}\,\mathrm{d}x\\
 &\leq& C_M^2\frac1{(1+2^j d(\xi, \xi'))^M}
|\bS^{d-1}| \int_0^1 r^{d-2}\,\mathrm{d}r.
\end{eqnarray*}
So
%
\begin{equation}
\forall M  \qquad| \langle\psi_{j,\xi}, \psi_{j,\xi'} \rangle
|\leq C'_M \frac1{(1+2^j d(\xi, \xi'))^M}.
\end{equation}

Now, let us choose $A_j$ as a maximal $K\eta$ net in the set $\chi
_j\cap\{x \in\bR^d,  |x|\le\frac1{\sqrt{2}}\}$ (case~1) or as a
maximal $K\eta$ net in the set $\chi_j\cap\{(1-2^{-2j}\le|x|\le1\}
$ (case 2).
Recall that $\eta\sim\frac\pi{2 2^{j+1}}$ and $K$ will be chosen later.

As, in case 1,
\[
\Card\{ \xi',   d (\xi' , \xi)
\sim Kl 2^{-j}\} \lesssim (Kl )^{d},\vspace*{-6pt}
\]
\[
\sum_{\xi' \neq\xi, \xi' \in A_j}
| \langle\psi_{j,\xi}, \psi_{j,\xi'} \rangle| \leq\sum
_{l=1}^{\fraca{2^j}K} (Kl )^{d}
C_M \frac1{(1+ Kl)^M} \leq C_M \sum_{l=1}^{\fraca{ 2^j}K} (Kl )^{d}
\frac1{ (Kl)^M} \leq\frac{2C_M}{K^{M-d}} \leq\alpha
\]
if $M-d \geq2$
and $K$ is large enough.
In case 2, again
\[
\Card\{ \xi',   d (\xi' , \xi)
\sim Kl 2^{-j}\} \lesssim (Kl )^{d-1}
\]
so
\begin{eqnarray*}
\sum_{\xi' \neq\xi, \xi' \in A_j}
| \langle\psi_{j,\xi}, \psi_{j,\xi'} \rangle| &\leq&\sum
_{l=1}^{\fraca{2^j}K} (Kl )^{d-1}
C_M \frac1{(1+ Kl)^M}\\
 &\leq& C_M \sum_{l=1}^{\fraca{ 2^j}K} (Kl )^{d-1}
\frac1{ (Kl)^M} \leq\frac{2C_M}{K^{M-d+1}} \leq\alpha
\end{eqnarray*}
if $M-d \geq1$
and $K$ is large enough.

 Hence, $ \cM(A_j)$ is invertible in both cases and we have
\[
c^{-1} \frac1{1-\alpha} \biggl( \sum_{\xi\in A_j} |\lambda_\xi
|^p \biggr)^{1/p} \leq
\biggl (\sum_{\xi'\in A_j} \biggl|\biggl \langle\sum_{\xi\in A_j} \lambda_\xi
 \langle\psi_{j,\xi} , \psi_{j,\xi'} \biggr\rangle\biggr|^p  \biggr)^{1/p}
\]
and
\[
 \biggl(\sum_{\xi\in A_j} | \lambda_\xi|^p \| \psi_{j,\xi} \|
^p_p \biggr)^{1/p} \lesssim \biggl\| \sum_{\xi\in A_j} \lambda_\xi\psi
_{j,\xi} \biggr\|_p \lesssim
\biggl (\sum_{\xi\in A_j} | \lambda_\xi|^p \| \psi_{j,\xi} \|
^p_p \biggr)^{1/p} .
\]

\subsection{Lower bounds associated sparse/dense cases and different
choices of~$A_j$}

Let $j$ be fixed and choose
\[
f =\sum_{\xi\in A_j} \beta_{j,\xi} \psi_{j,\xi}.
\]
We have
\[
f = \sum_{2^{j-1} <k <2^{j+1}} P_k(f),
\]
where $P_k$ is the orthogonal projector on $\cV_k(B^d)$.
So
\begin{eqnarray*}
\|R(f)\|^2&=& \langle R^*R(f),f \rangle= \sum_{2^{j-1} <k <2^{j+1}}
\langle\lambda^2_k P_k (f),f \rangle
\\
&\leq&  \biggl(\sup_{2^{j-1} <k <2^{j+1}} \lambda_k^2 \biggr) \sum_k \|
P_k(f) \|^2 \leq C 2^{-j(d-1)} \sum_{\xi\in A_j}
|\beta_{j,\xi}|^2.
\end{eqnarray*}

\subsubsection{Sparse choice, case 1}\label{sec8.2.1}

Let $f_i = \gamma\eps_i \psi_{j,\xi_i}$, $i\in\{1 , 2\}$, $\eps
_i$ is $+1$ or $-1$, in such a way
that
\[
\| f_1- f_2\|_p = \| \gamma\psi_{j,\xi_1}- \gamma\psi_{j,\xi_2}
\|_p = \gamma\| \psi_{j,\xi_1}- \psi_{j,\xi_2} \|_p \sim
\gamma(\| \psi_{j,\xi_1}\| +\|\psi_{j,\xi_2} \|_p).
\]
In case 1, $ \| \psi_{j,\xi}\|_r \sim2^{jd(1/2-1/r)}$.
So
\begin{eqnarray*}
f_i &\in& B^s_{\pi,r} (1)  \quad \Longleftrightarrow \quad  \gamma2^{jd(1/2-1/\pi)}
\sim2^{-js}
 \quad \Longleftrightarrow \quad  \gamma\sim2^{-j(s +d(1/2-1/\pi) ) },
\\
\delta&=&\| f_1- f_2\|_p \sim\gamma2^{jd(1/2-1/p)} \sim 2^{-j(s
+d(1/2-1/\pi)-d(1/2-1/p) }
=2^{-j(s - d(1/\pi-1/p) )}.
\end{eqnarray*}
On the other hand,
\[
K(P_1,P_2) = \frac12 \frac1{\epsilon^2}2^{-j(d-1)} \gamma^2 \sim
\frac12 \frac1{\epsilon^2}2^{-j(d-1)} 2^{-2j(s +d(1/2-1/\pi) )} =
\frac12 \frac1{\epsilon^2}
2^{-2j(s+d-1/2-d/\pi)}.
\]
Now, by the Fano inequality, if $j$ is chosen so that $\epsilon\sim
2^{-j(s+d-1/2-d/\pi)}$
(under the constraint
$s>d (1/\pi -(1-1/2d))$),
\[
\biggl (\frac2{\delta} \biggr)^p\bE\|f^\star-f_i\|_p^p \geq P(\|
f^\star-f_i\|_p >
\delta/2) \geq c.
\]
So, necessarily,
\[
\bE\|f^\star-f_i\|_p^p \geq c\delta^p \sim\epsilon^{\frace{s -
d(1/\pi-1/p) }{s+d-1/2-d/\pi}}.
\]

\begin{rem}
If
\[
d\biggl( 1/\pi- \biggl(1- \frac1{2d}\biggr)\biggr) < s \leq d(1/\pi-1/p)
\]
(so, necessarily, $\frac1p \leq1- \frac1{2d}),$
then
\[
\lim_{\epsilon\rightarrow0} \omega_p(s,\pi,q, M,\epsilon) \geq C>0.
\]
\end{rem}

\subsubsection{Sparse choice, case 2}\label{sec8.2.2}

In case 2,
$ \| \psi_{j,\xi}\|_r \sim2^{j (d+1)(1/2-1/r)}$,
so
\begin{eqnarray*}
f_i &\in& B^s_{\pi,r}(1)  \quad \Longleftrightarrow \quad  \gamma2^{j
(d+1)(1/2-1/\pi)} \sim2^{-js}
 \quad \Longleftrightarrow \quad  \gamma\sim2^{-j(s +(d+1)(1/2-1/\pi) )
 },
\\
\delta&=&\| f_1- f_2\|_p \sim\gamma2^{j (d+1)(1/2-1/p)} \sim 2^{-j(s
+(d+1)(1/2-1/\pi)-(d+1)(1/2-1/p) }\\
&=&2^{-j(s - (d+1)(1/\pi-1/p) )} .
\end{eqnarray*}
On the other hand,
\begin{eqnarray*}
K(P_1,P_2) &=& \frac12 \frac1{\epsilon^2}2^{-j(d-1)} \gamma^2 =
\frac12 \frac1{\epsilon^2}2^{-j(d-1)} 2^{-2j(s +(d+1)(1/2-1/\pi) )}\\
&\sim&\frac12 \frac1{\epsilon^2}
2^{-2j(s+d- (d+1)/\pi)}.
\end{eqnarray*}
Now, by the Fano inequality, if $\epsilon\sim2^{-j(s+d- (d+1)/\pi)}$
(under the constraint
$s> (d+1)(\frac1\pi-\frac d{d+1})$),
\[
 \biggl(\frac2{\delta} \biggr)^p\bE\|f^\star-f_i\|_p^p \geq P(\|
f^\star-f_i\|_p >
\delta/2) \geq c.
\]
So, necessarily,
\[
\bE\|f^\star-f_i\|_p^p \geq c\delta^p \sim\epsilon^{\fracc{(s -
(d+1)(1/\pi-1/p) )p}{s+d-(d+1)/\pi}}.
\]

\begin{rem}
If
\[
(d+1)\biggl(1/\pi-\frac d{d+1}\biggr) < s \leq (d+1)(1/\pi-1/p)
\]
(so, necessarily, $\frac1p < \frac d{d+1}$),
\[
\lim_{\epsilon\rightarrow0} \omega_p(s,\pi,q, M,\epsilon) \geq C>0.
\]
\end{rem}

\subsubsection{Dense choice, case 1}\label{sec8.2.3}

In this case, we take
\[
f_\rho= \gamma \sum_{\xi\in A_j} \epsilon_{\xi} \psi_{j,\xi},
\qquad\epsilon_\xi= \pm1,   \rho= (\epsilon_{\xi})_{\xi\in A_j}.
\]
As we are in case 1, we have
\[
\gamma^r\biggl\| \sum_{\xi\in A_j} \eps_\xi\psi_{j,\xi}\biggr\|_r^r
\sim \gamma^r 2^{jd( r/2-1)} \sum_{\xi\in A_j} | \eps_\xi
|^r\sim\gamma^r 2^{jd( r/2-1)} \Card(A_j)\sim\gamma^r 2^{jd r/2}.
\]
Using the Varshamov--Gilbert theorem (see \cite
{tsybakov08introdnonparestim}, Chapter 2), we consider a subset
$\cA$ of $\{-1,+1\}^{A_j}$ such that
$ \Card( \cA) \sim 2^{\fracd{1}{8} \Card( A_j)}$ and for $\rho\neq
\rho'$,  $\rho, \rho' \in\cA$,  $\| \rho-\rho' \|_1 \geq\frac
12 \Card( A_j)$.
Let us now restrict our set to
\begin{eqnarray*}
f_\rho&=& \gamma\sum_{\xi\in A_j} \epsilon_\xi
\psi_{j,\xi}, \qquad  \epsilon_\xi= \pm1,   \rho= (\epsilon_{\xi
})_{\xi
\in A_j},  \rho\in\cA,
\\
f_\rho&\in& B^{s}_{\pi,r}(1)  \quad \Longleftrightarrow \quad \gamma  \biggl(\sum
_{\xi\in A_j} \| \psi_{j,\xi}\|_\pi^\pi \biggr)^{1/\pi}\sim
\gamma2^{jd/2}
\sim2^{-js}.
\end{eqnarray*}
So we choose
\[
\gamma\sim2^{-j(s+d/2)}.
\]
Moreover,
\begin{eqnarray*}
\delta &=&\| f_\rho- f_{\rho'}\|_p= \gamma\biggl\| \sum_{\xi\in A_j}
(\epsilon_\xi-\epsilon'_\xi)\psi_{j,\xi} \biggr\|_p\sim\gamma  \biggl(
\sum_{\xi\in A_j} |\epsilon_\xi-\epsilon'_\xi
|^p\|\psi_{j,\xi} \|_p^p \biggr)^{1/p}
\\
&\sim&\gamma2^{jd(1/2-1/p)} \| \rho-\rho' \|_1^{1/p}
\sim2^{-j(s+d/2)} 2^{jd/2} =2^{-js}.
\end{eqnarray*}
Let us compute the Kullback distance,
\[
K(P_\rho, P_{\rho'}) = \frac1{2\epsilon^2} 2^{-j(d-1) } \| f_\rho-
f_{\rho'}\|_2^2
\sim
\frac1{2\epsilon^2} 2^{-j(d-1) } 2^{-2j(s+d/2)} 2^{jd} = \frac
1{2\epsilon^2}2^{-2j(s+d/2-1/2)},
\]
so, by the Fano inequality,
\[
\frac{\bE\|\hat f -f\|_p^p}{\delta^p} \geq 1/2 \wedge c 2^{\fracd{1}{8}
2^{jd}} \mathrm{e}^{-(\fracc1{2\epsilon^2})2^{-2j(s+d/2-1/2)}} \geq1/2
\]
if
\[
\epsilon\sim2^{-j(s+d-1/2)}.
\]
This implies
\[
\inf_{f \in B^{s}_{\pi,r}} \bE\|\hat f -f\|_p^p \geq c
\epsilon^{\fracc{sp}{s+d-1/2}}.\vspace*{-2pt}
\]

\subsubsection{Dense choice, case 2}\label{sec8.2.4}\vspace*{-2pt}

Similar to the previous case, we take now (with a slight abuse of
notation, since the subset $A$ obtained using the Varshamov--Gilbert
theorem is not the same $A$, as $A_j$ has also changed)
\[
f_\rho= \gamma\sum_{\xi\in A_j} \epsilon_\xi
\psi_{j,\xi}, \epsilon_\xi= \pm1,  \qquad   \rho= (\epsilon_{\xi
})_{\xi
\in A_j},  \rho\in\cA.
\]
As we are in case 2, we have
\begin{eqnarray*}
\gamma^r\biggl\| \sum_{\xi\in A_j} \eps_\xi\psi_{j,\xi}\biggr\|_r^r
&\sim& \gamma^r 2^{j(d+1)( r/2-1)} \sum_{\xi\in A_j} |\eps_\xi
|^r\\[-2pt]
&\sim&\gamma^r 2^{j(d+1)( r/2-1)} \Card(A_j)\sim\gamma^r
2^{j[(d+1) r/2-2]},
\\[-2pt]
f_\rho\in B^{s}_{\pi,r}(1)  \quad \Longleftrightarrow \quad \gamma  \biggl(\sum
_{\xi\in A_j} \| \psi_{j,\xi}\|_\pi^\pi \biggr)^{1/\pi} &\sim&
\gamma \bigl( 2^{j(d-1)} 2^{j(d+1) (\pi/2-1) }  \bigr)^{1/\pi}\\[-2pt]
&\sim& \gamma2^{-j(\fracb{d+1}2 -\fraca2\pi)}
\sim2^{-js}.
\end{eqnarray*}
So we choose
\[
\gamma\sim2^{-j(s+\fracb{d+1}2 -\fraca2\pi)}.
\]
 Moreover,
\begin{eqnarray*}
 \delta &=&\| f_\rho- f_{\rho'}\|_p= \gamma\biggl\| \sum_{\xi\in A_j}
(\epsilon_\xi-\epsilon'_\xi)\psi_{j,\xi} \biggr\|_p\sim\gamma  \biggl(
\sum_{\xi\in A_j} |\epsilon_\xi-\epsilon'_\xi
|^p\|\psi_{j,\xi} \|_p^p \biggr)^{1/p}\\[-2pt]
&\sim&\gamma2^{j(d+1)(1/2-1/p)} \| \rho-\rho' \|_1^{1/p}
\sim\gamma2^{j(d+1)(1/2-1/p)} 2^{j(d-1) \fraca1p}\\[-2pt]
& \sim& 2^{-j(s+\fracb{d+1}2 -\fraca2\pi)} 2^{j(\fracb{d+1}2 -\fraca
2p)}=2^{-j(s -2(1/\pi-1/p) )}.
\end{eqnarray*}

Let us compute the Kullback distance:
\begin{eqnarray*}
K(P_\rho, P_{\rho'}) &=& \frac1{2\epsilon^2} 2^{-j(d-1) } \| f_\rho-
f_{\rho'}\|_2^2\\[-2pt]
&\sim&
\frac1{2\epsilon^2} 2^{-j(d-1) } 2^{-2j(s+\fracb{d+1}2 -\fraca2\pi)}
2^{j(d-1)} = \frac1{2\epsilon^2} 2^{-2j(s+\fracb{d+1}2 -\fraca2\pi),}
\end{eqnarray*}
so, by the Fano inequality,
\[
\frac{\bE\|\hat f -f\|_p^p}{\delta^p} \geq 1/2 \wedge c 2^{\fracd
{1}{8} 2^{j(d-1)}} \mathrm{e}^{-\fracc1{2\epsilon^2}
2^{-2j(s+\fracb{d+1}2 -\fraca2\pi )}}
\geq1/2
\]
if
\[
\epsilon\sim2^{-j(s+d-2/\pi)}.
\]
This implies
\[
\inf_{f \in B^{s}_{\pi,r}(1)} \bE\|\hat f -f\|_p^p \geq c
\epsilon^{\fracc{p (s-2(1/\pi-1/p) )}{s+d- 2/\pi}}.
\]

\begin{rem}The case $p=\infty$ can be handled using the same arguments
without difficulties.
\end{rem}


\printhistory


\begin{thebibliography}{29}

\bibitem{AAR}
\begin{bbook}[mr]
\bauthor{\bsnm{Andrews},~\bfnm{George~E.}\binits{G.E.}},
  \bauthor{\bsnm{Askey},~\bfnm{Richard}\binits{R.}} \AND
  \bauthor{\bsnm{Roy},~\bfnm{Ranjan}\binits{R.}}
(\byear{1999}).
\btitle{Special Functions}.
\bseries{Encyclopedia of Mathematics and Its Applications}
\bvolume{71}.
\baddress{Cambridge}: \bpublisher{Cambridge Univ. Press}.
\bid{mr={1688958}}
\bptnote{check year}
\end{bbook}
\endbibitem

\bibitem{baldi-2008}
\begin{barticle}[mr]
\bauthor{\bsnm{Baldi},~\bfnm{P.}\binits{P.}},
  \bauthor{\bsnm{Kerkyacharian},~\bfnm{G.}\binits{G.}},
  \bauthor{\bsnm{Marinucci},~\bfnm{D.}\binits{D.}} \AND
  \bauthor{\bsnm{Picard},~\bfnm{D.}\binits{D.}}
(\byear{2009}).
\btitle{Adaptive density estimation for directional data using needlets}.
\bjournal{Ann. Statist.}
\bvolume{37}
\bpages{3362--3395}.
\bid{doi={10.1214/09-AOS682}, issn={0090-5364}, mr={2549563}}
\end{barticle}
\endbibitem

\bibitem{baldi09asymp}
\begin{barticle}[mr]
\bauthor{\bsnm{Baldi},~\bfnm{P.}\binits{P.}},
  \bauthor{\bsnm{Kerkyacharian},~\bfnm{G.}\binits{G.}},
  \bauthor{\bsnm{Marinucci},~\bfnm{D.}\binits{D.}} \AND
  \bauthor{\bsnm{Picard},~\bfnm{D.}\binits{D.}}
(\byear{2009}).
\btitle{Asymptotics for spherical needlets}.
\bjournal{Ann. Statist.}
\bvolume{37}
\bpages{1150--1171}.
\bid{doi={10.1214/08-AOS601}, issn={0090-5364}, mr={2509070}}
\end{barticle}
\endbibitem

\bibitem{birge}
\begin{bmisc}[auto:STB|2011-03-03|12:04:44]
\bauthor{\bsnm{Birge},~\bfnm{L.}\binits{L.}}
(\byear{2001}).
\bhowpublished{A new look at an old result: Fano's lemma. Prepublication 632,
  LPMA}.
\end{bmisc}
\endbibitem

\bibitem{Davison}
\begin{barticle}[mr]
\bauthor{\bsnm{Davison},~\bfnm{M.~E.}\binits{M.E.}}
(\byear{1981}).
\btitle{A singular value decomposition for the {R}adon transform in
  {$n$}-dimensional {E}uclidean space}.
\bjournal{Numer. Funct. Anal. Optim.}
\bvolume{3}
\bpages{321--340}.
\bid{doi={10.1080/01630568108816093}, issn={0163-0563}, mr={0629949}}
\end{barticle}
\endbibitem

\bibitem{DUXU}
\begin{bbook}[mr]
\bauthor{\bsnm{Dunkl},~\bfnm{Charles~F.}\binits{C.F.}} \AND
  \bauthor{\bsnm{Xu},~\bfnm{Yuan}\binits{Y.}}
(\byear{2001}).
\btitle{Orthogonal Polynomials of Several Variables}.
\bseries{Encyclopedia of Mathematics and Its Applications}
\bvolume{81}.
\baddress{Cambridge}: \bpublisher{Cambridge Univ. Press}.
\bid{doi={10.1017/CBO9780511565717}, mr={1827871}}
\end{bbook}
\endbibitem

\bibitem{Erdelyi}
\begin{bbook}[mr]
\bauthor{\bsnm{Erd{\'e}lyi},~\bfnm{Arthur}\binits{A.}},
  \bauthor{\bsnm{Magnus},~\bfnm{Wilhelm}\binits{W.}},
  \bauthor{\bsnm{Oberhettinger},~\bfnm{Fritz}\binits{F.}} \AND
  \bauthor{\bsnm{Tricomi},~\bfnm{Francesco~G.}\binits{F.G.}}
(\byear{1981}).
\btitle{Higher Transcendental Functions. {V}ol. {II}}.
\baddress{Melbourne, FL}: \bpublisher{Robert E. Krieger Publishing Co. Inc.}
\bid{mr={0698780}}
\end{bbook}
\endbibitem

\bibitem{fano}
\begin{bmisc}[auto:STB|2011-03-03|12:04:44]
\bauthor{\bsnm{Fano},~\bfnm{R.}\binits{R.}}
(\byear{1952}).
\bhowpublished{Class notes for transmission of information, course. 6.574.
 MIT,  Cambridge, MA}.
\end{bmisc}
\endbibitem

\bibitem{grafakos}
\begin{bbook}[mr]
\bauthor{\bsnm{Grafakos},~\bfnm{Loukas}\binits{L.}}
(\byear{2004}).
\btitle{Classical and Modern {F}ourier Analysis}.
\baddress{Upper Saddle River, NJ}: \bpublisher{Pearson Education, Inc}.
\bid{mr={2449250}}
\end{bbook}
\endbibitem

\bibitem{Helgason}
\begin{bbook}[mr]
\bauthor{\bsnm{Helgason},~\bfnm{Sigurdur}\binits{S.}}
(\byear{1999}).
\btitle{The {R}adon Transform},
\bedition{2nd} ed.
\bseries{Progress in Mathematics}
\bvolume{5}.
\baddress{Boston, MA}: \bpublisher{Birkh\"auser}.
\bid{mr={1723736}}
\end{bbook}
\endbibitem

\bibitem{5authors}
\begin{barticle}[mr]
\bauthor{\bsnm{Kerkyacharian},~\bfnm{G{\'e}rard}\binits{G.}},
  \bauthor{\bsnm{Kyriazis},~\bfnm{George}\binits{G.}},
  \bauthor{\bsnm{Le~Pennec},~\bfnm{Erwan}\binits{E.}},
  \bauthor{\bsnm{Petrushev},~\bfnm{Pencho}\binits{P.}} \AND
  \bauthor{\bsnm{Picard},~\bfnm{Dominique}\binits{D.}}
(\byear{2010}).
\btitle{Inversion of noisy {R}adon transform by {SVD} based needlets}.
\bjournal{Appl. Comput. Harmon. Anal.}
\bvolume{28}
\bpages{24--45}.
\bid{doi={10.1016/j.acha.2009.06.001}, issn={1063-5203}, mr={2563258}}
\end{barticle}
\endbibitem

\bibitem{kyoto}
\begin{bincollection}[auto:STB|2011-03-03|12:04:44]
\bauthor{\bsnm{Kerkyacharian},~\bfnm{G.}\binits{G.}} \AND
  \bauthor{\bsnm{Picard},~\bfnm{D.}\binits{D.}}
(\byear{2009}).
\btitle{New generation wavelets associated with statistical problems}.
In \bbooktitle{The 8th Workshop on Stochastic Numerics}
\bpages{119--146}.
\baddress{Kyoto Univ.}: \bpublisher{Research Institute for Mathematical
  Sciences}.
\end{bincollection}
\endbibitem

\bibitem{pxukball}
\begin{barticle}[mr]
\bauthor{\bsnm{Kyriazis},~\bfnm{G.}\binits{G.}},
  \bauthor{\bsnm{Petrushev},~\bfnm{P.}\binits{P.}} \AND
  \bauthor{\bsnm{Xu},~\bfnm{Yuan}\binits{Y.}}
(\byear{2008}).
\btitle{Decomposition of weighted {T}riebel--{L}izorkin and {B}esov spaces on
  the ball}.
\bjournal{Proc. Lond. Math. Soc. (3)}
\bvolume{97}
\bpages{477--513}.
\bid{doi={10.1112/plms/pdn010}, issn={0024-6115}, mr={2439670}}
\end{barticle}
\endbibitem

\bibitem{LOG}
\begin{barticle}[mr]
\bauthor{\bsnm{Logan},~\bfnm{B.~F.}\binits{B.F.}} \AND
  \bauthor{\bsnm{Shepp},~\bfnm{L.~A.}\binits{L.A.}}
(\byear{1975}).
\btitle{Optimal reconstruction of a function from its projections}.
\bjournal{Duke Math. J.}
\bvolume{42}
\bpages{645--659}.
\bid{issn={0012-7094}, mr={0397240}}
\end{barticle}
\endbibitem

\bibitem{Louis}
\begin{barticle}[mr]
\bauthor{\bsnm{Louis},~\bfnm{Alfred~K.}\binits{A.K.}}
(\byear{1984}).
\btitle{Orthogonal function series expansions and the null space of the {R}adon
  transform}.
\bjournal{SIAM J. Math. Anal.}
\bvolume{15}
\bpages{621--633}.
\bid{doi={10.1137/0515047}, issn={0036-1410}, mr={0740700}}
\end{barticle}
\endbibitem

\bibitem{muciaccia97fastspherharmonanaly}
\begin{barticle}[auto:STB|2011-03-03|12:04:44]
\bauthor{\bsnm{Muciaccia},~\bfnm{P.~F.}\binits{P.F.}},
  \bauthor{\bsnm{Natoli},~\bfnm{P.}\binits{P.}} \AND
  \bauthor{\bsnm{Vittorio},~\bfnm{N.}\binits{N.}}
(\byear{1997}).
\btitle{Fast spherical harmonic analysis: A quick algorithm for generating
  and/or inverting full-sky high-resolution cosmic microwave background
  anisotropy maps}.
\bjournal{The Astrophysical J. Lett.}
\bvolume{488}
\bpages{63--66}.
\end{barticle}
\endbibitem

\bibitem{pnarco}
\begin{barticle}[mr]
\bauthor{\bsnm{Narcowich},~\bfnm{F.}\binits{F.}},
  \bauthor{\bsnm{Petrushev},~\bfnm{P.}\binits{P.}} \AND
  \bauthor{\bsnm{Ward},~\bfnm{J.}\binits{J.}}
(\byear{2006}).
\btitle{Decomposition of {B}esov and {T}riebel--{L}izorkin spaces on the
  sphere}.
\bjournal{J. Funct. Anal.}
\bvolume{238}
\bpages{530--564}.
\bid{issn={0022-1236}, mr={2253732}}
\end{barticle}
\endbibitem

\bibitem{NPW}
\begin{barticle}[mr]
\bauthor{\bsnm{Narcowich},~\bfnm{F.~J.}\binits{F.J.}},
  \bauthor{\bsnm{Petrushev},~\bfnm{P.}\binits{P.}} \AND
  \bauthor{\bsnm{Ward},~\bfnm{J.~D.}\binits{J.D.}}
(\byear{2006}).
\btitle{Localized tight frames on spheres}.
\bjournal{SIAM J. Math. Anal.}
\bvolume{38}
\bpages{574--594 (electronic)}.
\bid{doi={10.1137/040614359}, issn={0036-1410}, mr={2237162}}
\end{barticle}
\endbibitem

\bibitem{NATT}
\begin{bbook}[mr]
\bauthor{\bsnm{Natterer},~\bfnm{F.}\binits{F.}}
(\byear{2001}).
\btitle{The Mathematics of Computerized Tomography}.
\bseries{Classics in Applied Mathematics}
\bvolume{32}.
\baddress{Philadelphia, PA}: \bpublisher{SIAM}.
\bnote{Reprint of the 1986 original}.
\bid{mr={1847845}}
\end{bbook}
\endbibitem

\bibitem{PXU}
\begin{barticle}[mr]
\bauthor{\bsnm{Petrushev},~\bfnm{Pencho}\binits{P.}} \AND
  \bauthor{\bsnm{Xu},~\bfnm{Yuan}\binits{Y.}}
(\byear{2005}).
\btitle{Localized polynomial frames on the interval with {J}acobi weights}.
\bjournal{J.~Fourier Anal. Appl.}
\bvolume{11}
\bpages{557--575}.
\bid{doi={10.1007/s00041-005-4072-3}, issn={1069-5869}, mr={2182635}}
\end{barticle}
\endbibitem

\bibitem{pxuball}
\begin{barticle}[mr]
\bauthor{\bsnm{Petrushev},~\bfnm{Pencho}\binits{P.}} \AND
  \bauthor{\bsnm{Xu},~\bfnm{Yuan}\binits{Y.}}
(\byear{2008}).
\btitle{Localized polynomial frames on the ball}.
\bjournal{Constr. Approx.}
\bvolume{27}
\bpages{121--148}.
\bid{doi={10.1007/s00365-007-0678-9}, issn={0176-4276}, mr={2336420}}
\end{barticle}
\endbibitem

\bibitem{petrush}
\begin{barticle}[mr]
\bauthor{\bsnm{Petrushev},~\bfnm{Pencho~P.}\binits{P.P.}}
(\byear{1999}).
\btitle{Approximation by ridge functions and neural networks}.
\bjournal{SIAM J. Math. Anal.}
\bvolume{30}
\bpages{155--189 (electronic)}.
\bid{doi={10.1137/S0036141097322959}, issn={0036-1410}, mr={1646689}}
\end{barticle}
\endbibitem

\bibitem{Pisierentropie}
\begin{bincollection}[mr]
\bauthor{\bsnm{Pisier},~\bfnm{Gilles}\binits{G.}}
(\byear{1983}).
\btitle{Some applications of the metric entropy condition to harmonic
  analysis}.
In \bbooktitle{Banach Spaces, Harmonic Analysis, and Probability Theory
  ({S}torrs, {C}onn., 1980/1981)}.
\bseries{Lecture Notes in Math.}
\bvolume{995}
\bpages{123--154}.
\baddress{Berlin}: \bpublisher{Springer}.
\bid{mr={0717231}}
\end{bincollection}
\endbibitem

\bibitem{STW}
\begin{bbook}[mr]
\bauthor{\bsnm{Stein},~\bfnm{Elias~M.}\binits{E.M.}} \AND
  \bauthor{\bsnm{Weiss},~\bfnm{Guido}\binits{G.}}
(\byear{1971}).
\btitle{Introduction to {F}ourier Analysis on {E}uclidean Spaces}.
\bseries{Princeton Mathematical Series}
\bvolume{32}.
\baddress{Princeton, NJ}: \bpublisher{Princeton Univ. Press}.
\bid{mr={0304972}}
\end{bbook}
\endbibitem

\bibitem{SZG}
\begin{bbook}[mr]
\bauthor{\bsnm{Szeg{\H{o}}},~\bfnm{G{\'a}bor}\binits{G.}}
(\byear{1975}).
\btitle{Orthogonal Polynomials}, \bedition{4th} ed.
\bseries{American Mathematical Society, Colloquium Publications}
\bvolume{XXIII}.
\baddress{Providence, RI}: \bpublisher{Amer. Math. Soc.}
\bid{mr={0372517}}
\end{bbook}
\endbibitem

\bibitem{tsybakov08introdnonparestim}
\begin{bbook}[auto:STB|2011-03-03|12:04:44]
\bauthor{\bsnm{Tsybakov},~\bfnm{Alexandre~B.}\binits{A.B.}}
(\byear{2008}).
\btitle{Introduction to Nonparametric Estimation}.
\baddress{Berlin}: \bpublisher{Springer}.
\end{bbook}
\endbibitem

\bibitem{cubxu}
\begin{barticle}[mr]
\bauthor{\bsnm{Xu},~\bfnm{Yuan}\binits{Y.}}
(\byear{1998}).
\btitle{Orthogonal polynomials and cubature formulae on spheres and on balls}.
\bjournal{SIAM J. Math. Anal.}
\bvolume{29}
\bpages{779--793 (electronic)}.
\bid{doi={10.1137/S0036141096307357}, issn={0036-1410}, mr={1617720}}
\end{barticle}
\endbibitem

\bibitem{xu7}
\begin{barticle}[mr]
\bauthor{\bsnm{Xu},~\bfnm{Yuan}\binits{Y.}}
(\byear{2007}).
\btitle{Reconstruction from {R}adon projections and orthogonal expansion on a
  ball}.
\bjournal{J. Phys. A}
\bvolume{40}
\bpages{7239--7253}.
\bid{doi={10.1088/1751-8113/40/26/010}, issn={1751-8113}, mr={2344454}}
\end{barticle}
\endbibitem

\end{thebibliography}
\end{document}